\documentclass{article}

\usepackage{arxiv}

\usepackage[utf8]{inputenc} 
\usepackage[T1]{fontenc}    
\usepackage{hyperref}       
\usepackage{url}            
\usepackage{booktabs}       
\usepackage{amsfonts}       
\usepackage{nicefrac}       
\usepackage{microtype}      
\usepackage{lipsum}
\usepackage{graphicx}
\graphicspath{ {./images/} }

\usepackage{graphics}
\usepackage{amsmath} 
\usepackage{float} 

\usepackage[title,titletoc]{appendix}
\usepackage{latexsym}
\usepackage{subfigure}
\allowdisplaybreaks[4]
\usepackage{algorithm}
\usepackage{algpseudocode}
\usepackage{epsfig}
\usepackage{color}
\usepackage{amssymb}
\usepackage{ntheorem}

\newtheorem{remark}{Remark}[section]
\usepackage{booktabs} 
\usepackage{fancyhdr} 
\usepackage{bm}
\usepackage{tabularx}
\usepackage[textsize=tiny,tickmarkheight=5pt]{todonotes}
\numberwithin{equation}{section} 
\numberwithin{figure}{section}

\newcommand{\half}{\frac{1}{2}}
\newcommand{\jmh}{{j-\half}}
\newcommand{\jph}{{j+\half}}
\newcommand{\kmh}{{k-\half}}
\newcommand{\kph}{{k+\half}}
\newcommand{\Cjk}{C_{j,k}}

\title{A low-dissipation central scheme for ideal MHD}

\author{
 Yu-Chen Cheng \\
  Institute of Mathematics\\ 
  University of Würzburg\\ 
  Würzburg, 97074\\
  Germany\\
  \texttt{yu-chen.cheng@stud-mail.uni-wuerzburg.de} \\
   \And
 Praveen Chandrashekar \\
  Centre for Applicable Mathematics\\
  Tata Institute of Fundamental Research \\
  Bengaluru, 560065\\
  India \\
  \texttt{praveen@tifrbng.res.in} \\
  \And
 Christian Klingenberg \\
  Institute of Mathematics\\ 
  University of Würzburg\\ 
  Würzburg, 97074\\
  Germany\\
  \texttt{christian.klingenberg@uni-wuerzburg.de} \\
}

\begin{document}

\maketitle

\begin{abstract}
Central schemes for conservation laws are Riemann solver free methods which are simple and easy to implement. In recent work for Euler equations~\cite{LDCU1} their accuracy has been enhanced in terms of better resolution of contact waves. In this paper, we extend this low dissipation central upwind method to the ideal MHD system in one- and two-dimensions.  In the two-dimensional case, we separate the variables into two groups: hydrodynamic and magnetic, which are stored at cell centers and faces, respectively.  For the  the hydrodynamic variables, we apply the low dissipation central upwind scheme while for the magnetic variables, a constrained transport method is used which maintains the divergence-free property of the magnetic field. The time integration is performed with third order strong stability preserving Runge-Kutta scheme. To validate the proposed scheme, we apply this method to several challenging test cases. The results show that the LDCU correction term plays a useful role at the contact discontinuity and enhances the resolution of waves.  We also observe experimental second-order accuracy for smooth solutions and the divergence-free condition is maintained to machine precision.
\end{abstract}


\section{Introduction}
We consider the numerical solution of hyperbolic conservation laws arising in ideal magnetohydrodynamics in one and two dimensions by central type of schemes. The numerical solution of MHD equations is complicated by the fact that the divergence-free condition on the magnetic field $\textbf{B}$ has to be maintained in the discrete solution in order to have stable computations.  To maintain  the condition $\nabla\cdot \textbf{B}=0$ for all times, there are different methods available in the literature;

\begin{itemize}
    \item
    the Godunov-Powell method (or eight-wave formulation) in \cite{Powell} adds a non-conservative source term to recover the missing eigenvector and reduce the divergence errors. More references can be found in~\cite{Powell_intro1},~\cite{Powell_intro2}. This also facilitates the construction of entropy stable schemes as in~\cite{Powell_intro3},~\cite{Powell_intro4}.
    \item the Leray projection method \cite{Project_intro}, which solves for a potential $\phi$ from the Poisson equation $\Delta\phi=-\nabla\cdot \mathbf{B}$ and updates the magnetic field $\mathbf{B}$ obtained from a numerical scheme to the divergence-free magnetic field $\mathbf{B}^c$ by $\mathbf{B}^c = \mathbf{B}+\nabla\phi$ at the end of each time step. The scheme in~\cite{BT1} proposes a central scheme with this Leray projection method. 
    \item the divergence cleaning method in \cite{GLM} is an approach coupling the divergence constraint with the evolution equations using a generalized Lagrange multiplier (GLM).  The core of the method is introducing an auxiliary correction function $\psi$ to allow the divergence errors propagate to the domain boundaries and are damped at the same time. 
    \item the constrained transport method (see \cite{CT_1st}). The idea of this method is using the approximation of the electric field $\Omega$ at the staggered grids to preserve a specific discretization of divergence by the induction equation $\mathbf{B}_t = -\nabla\times\Omega$.     Under this framework, there are some variants \cite{CT_intro1},~\cite{CT_intro2},~\cite{CT_intro3},~\cite{CT_intro4},~\cite{CT_intro5},~\cite{Toth}. An approach to construct divergence-free method is introduced in~\cite{UCT} using the upwind-constrained transport (UCT) idea. The UCT-HLL scheme approximates the electric field at staggered grid and is systematically developed in~\cite{UCT-HLL1}. An alternative way of the UCT-HLL scheme is shown in~\cite{UCT-HLL2}.
\end{itemize}

In this work, we consider central schemes which are Riemann solver-free and combine it with a constrained transport method. In \cite{KT1}, the central scheme is designed to use the narrower interval over Riemann fans to approximate the average at the next time step, instead of estimating the interval over the same length as the computational cell. Later, \cite{KT2} modified the KT scheme that considers more specific information to define the range of the interval over Riemann fans. Based on the same framework, there are some varieties of central schemes in \cite{central_not_rectangle}, \cite{central_1}. For example, the scheme in \cite{central_differ_minmod} changes the interpolation approach to achieve higher resolution.   In~\cite{LDCU1}, a low-dissipation central-upwind scheme (LDCU) is further developed to reduce dissipation at contact discontinuities using the property of contact waves.  The second version of this scheme in~\cite{LDCU2} considers a more likely position where the contact occurs in order to decrease the oscillation at the boundary. In this paper we generalize the LDCU method to MHD in a conservative way while guaranteeing the constraint of the solenoidal magnetic field.  We combine the central scheme with a constrained transport method.  

The rest of the paper is organized as follows.  In Section 2, we first summarize the fully discrete LDCU scheme~\cite{LDCU2} for the 1-D MHD equations, and then derive the 1-D semi-discrete scheme. Next, with the help of dimension-by-dimension method and the constrained transport method in~\cite{UCT-HLL2}, we construct a 2-D semi-discrete scheme in Section 3.  In Section 4, we implement our 1-D fully-discrete scheme and 2-D semi-discrete scheme to several numerical experiments, and demonstrate that the results are consistent with our claim. Finally, we end up with summary and conclusions in Section 5.

\section{1-D LDCU MHD scheme}
 
In this section, we extend the basic idea of the LDCU scheme for the 1-D Euler equations proposed by~\cite{LDCU1} and~\cite{LDCU2} to the MHD model.  Following the structure of the KT-type central scheme, e.g.~\cite{KT1},~\cite{KT2}, or~\cite{LDCU1}, we separate the scheme into three stages: reconstruction, evolution, and projection.  The most important difference between the LDCU scheme and the other central schemes occurs in the projection step.  The projection step uses more information of the contact wave to simulate a realistic jump at discontinuity, when projecting the averaged value from the staggered grids back to the uniform grids.  In the MHD model, we need to consider more quantities and their properties at contact wave.  

\subsection{1-D LDCU scheme for the MHD system} \label{sec.1DLDCU}
\begin{figure}[htbp]
    \centering
    \subfigure[Step 1: Reconstruction at $t^n$]{
    \includegraphics[width=0.51\linewidth]{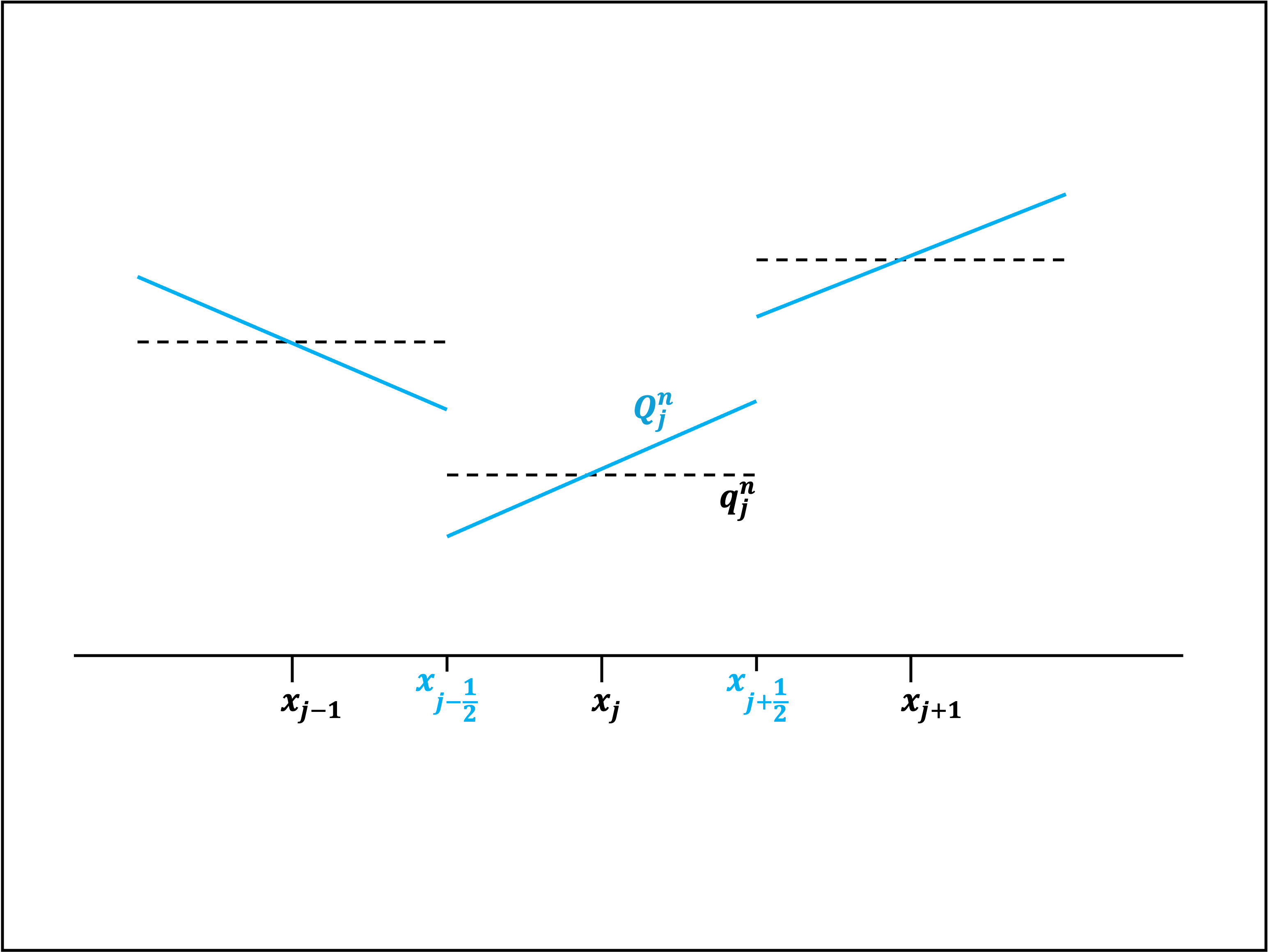}
    }\\
    \subfigure[Step 2: Evolution from $t^n$ to $t^{n+1}$]{
    \includegraphics[width=0.51\linewidth]{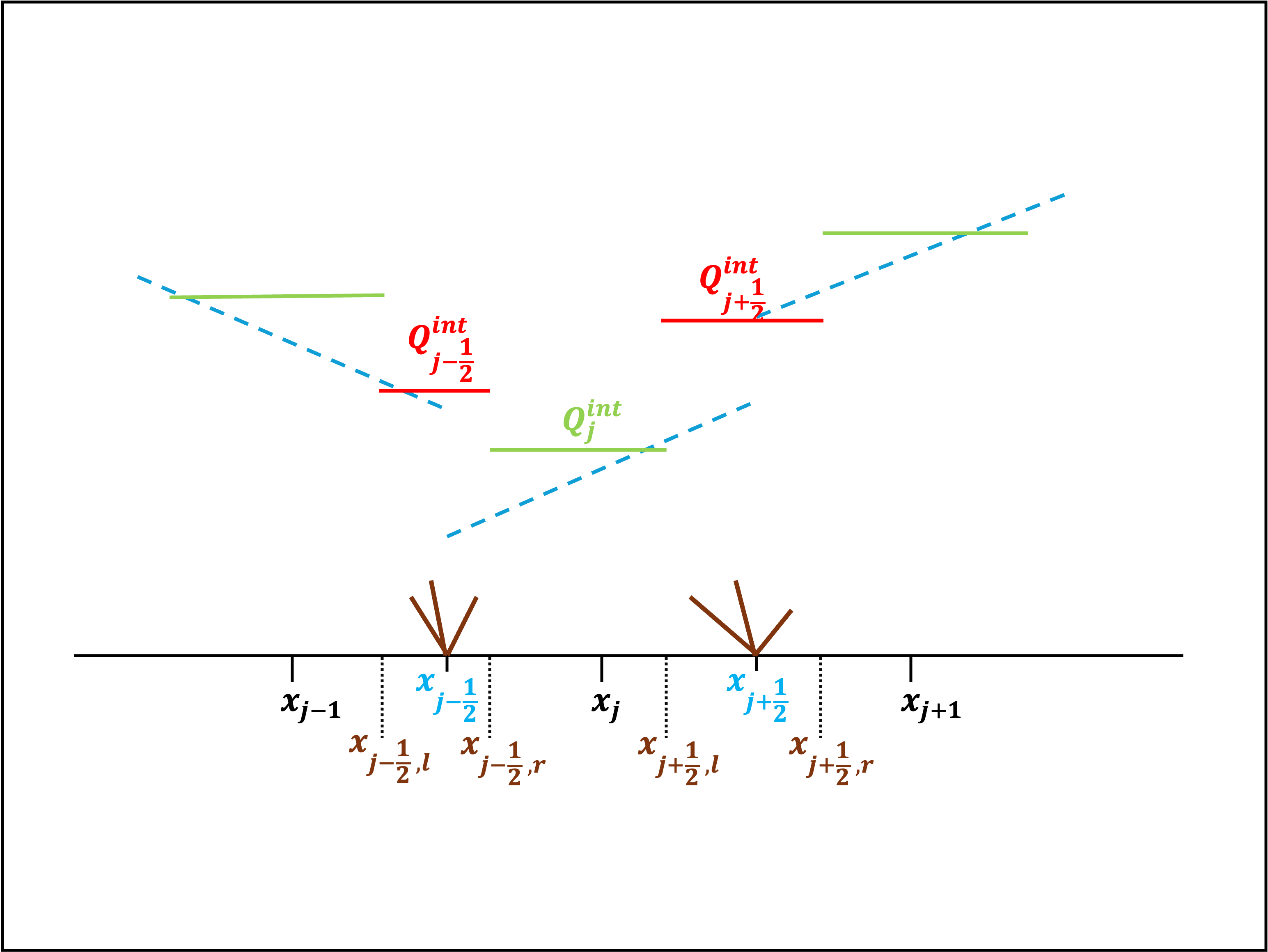}
    }\\
    \subfigure[Step 3: Projection at $t^{n+1}$]{
    \includegraphics[width=0.51\linewidth]{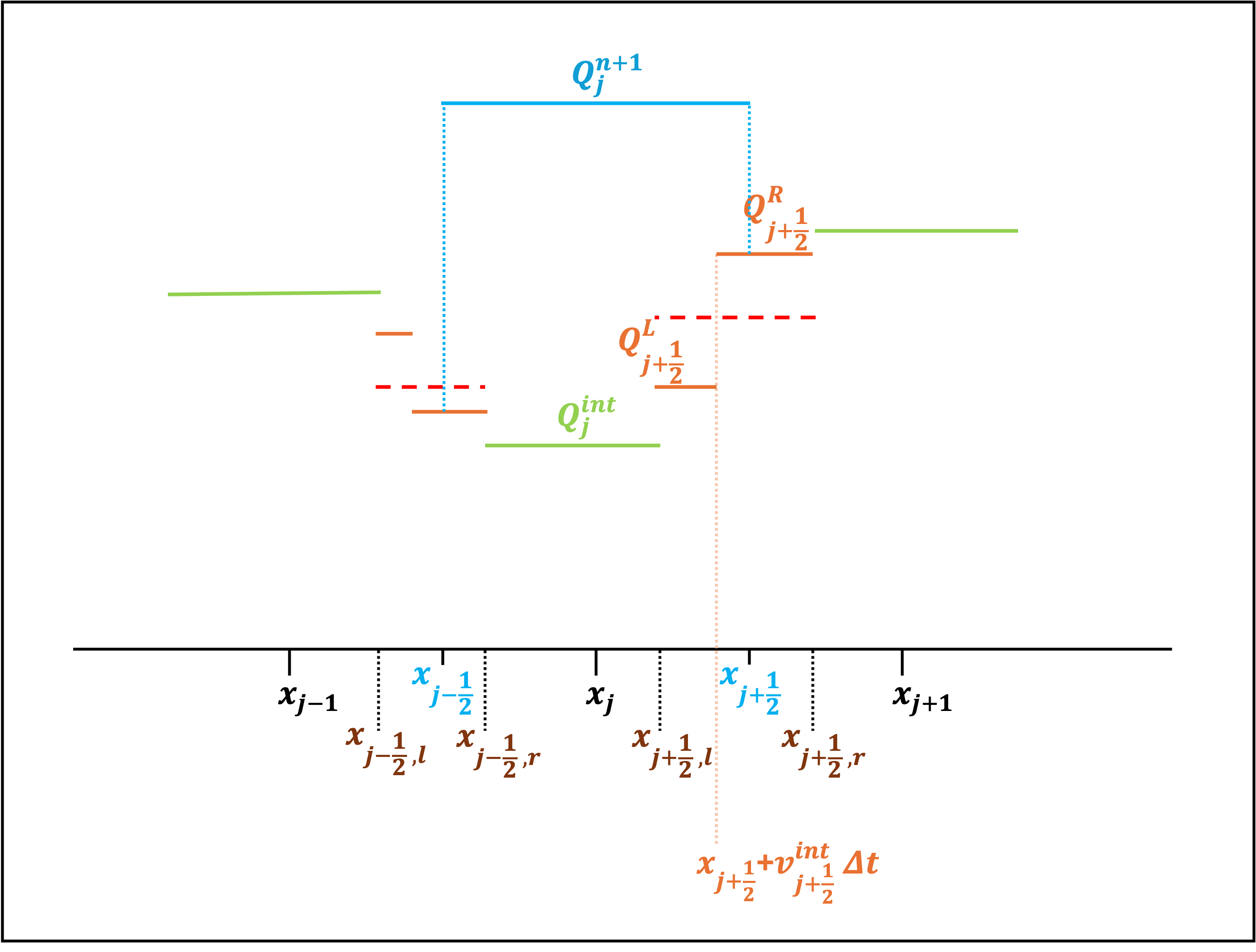}
    }
    \caption{Three steps of the LDCU method, following~\cite{LDCU2}.}
    \label{fig.LDCU step}
\end{figure}

In the MHD system, in addition to the hydrodynamic variables density $\rho$, momentum $(\rho v_i)$, $i=1,2,3$, and total energy $E$, we also have the magnetic variables $B_i$ to which the LDCU scheme is applied.   Consider 1-D MHD system
\begin{equation} \label{1Dlaw}
q_t + f(q)_x = 0,
\end{equation}
where
\begin{equation*}
    q = \left[
    \begin{matrix}
        \rho  \\
        \rho v_1 \\
        \rho v_2 \\
        \rho v_3 \\
        B_2 \\
        B_3 \\
        E
    \end{matrix}
    \right],\qquad
    f(q) = \left[
    \begin{matrix}
        \rho v_1  \\
        \rho v_1^2 + p + \frac{1}{2}|\textbf{B}|^2 - B_1^2 \\
        \rho v_1v_2 - B_1B_2 \\
        \rho v_1v_3 - B_1B_3 \\
        B_2v_1 - B_1v_2 \\
        B_3v_1 - B_1v_3 \\
        (E+p+\frac{1}{2}|\textbf{B}|^2)v_1 - B_1(B_1v_1+B_2v_2+B_3v_3)
    \end{matrix}
    \right].
\end{equation*}
with $\textbf{B}=(B_1, B_2, B_3)^T$ and $B_1 = $ constant.  Consider a partitioning of the domain into disjoint cells $C_j = (x_\jmh, x_\jph)$ with $\Delta x = x_\jph - x_\jmh$ and $x_j = (x_\jmh + x_\jph)/2$. At the time $t^n$ the solution is made of piecewise constant states given by the cell averages $Q_j^n$.  Similarly to the 1-D Euler model, we firstly consider a reconstruction at the first step and secondly approximate the intermediate average over the Riemann fans in the second step. In the third projection step, we use the the algorithm in~\cite{LDCU2} Section 2.2 to reconstruct the density at both side of the contact discontinuity.  Afterwards, based on the property of contact wave, we approximate the momentum and the energy by the reconstructed density. The framework of the scheme is shown in Figure~\ref{fig.LDCU step}, and we illustrate the details in several steps as follows. 

~\\
\noindent\textbf{Step 1: Reconstruction:}
In the first step, we build a piecewise linear, second-order reconstruction,
\begin{equation} \label{1DRec}
    Q_j(x,t^n) = Q^n_j + (x-x_j)(Q_x)^n_j, \qquad x \in C_j
\end{equation} 
where the slopes $(Q_x)_j^n$ are approximated using a limiter to ensure some non-oscillatory property. 

\begin{figure}[htbp]
    \centering
    \includegraphics[width=0.7\linewidth]{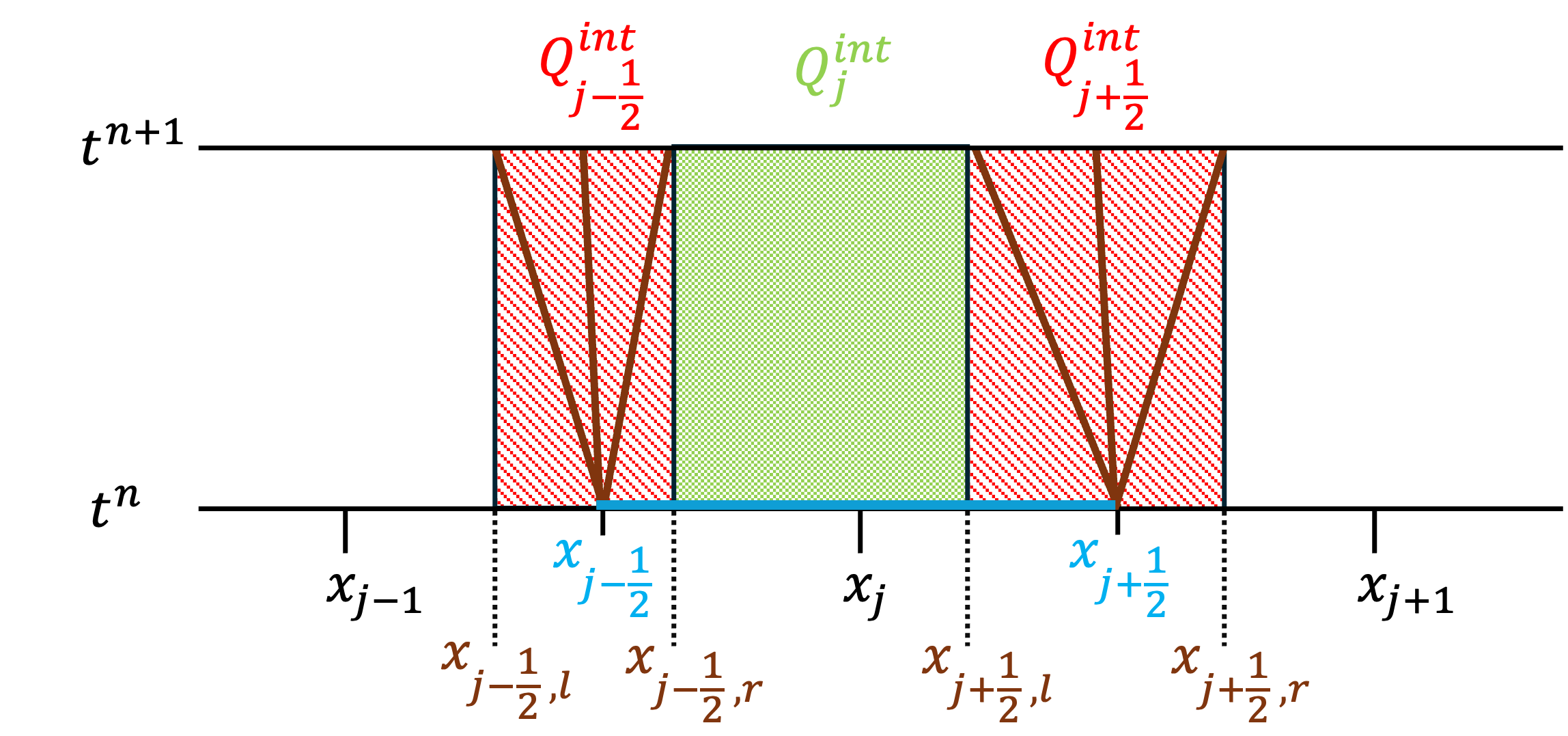}
    \caption{Smooth (green) and unsmooth (red) region. This figure is inspired by Figure 3 in \cite{LDCU1}.}
    \label{unORsmooth}
\end{figure}

~\\
\noindent\textbf{Step 2: Evolution:}
Evaluating the piecewise linear reconstructions at the cell interfaces gives the two states 
\begin{equation*}
    Q^{n,-}_{j+\frac{1}{2}} = Q^n_j + \frac{\Delta x}{2}(Q_x)^n_j,\qquad
    Q^{n,+}_{j+\frac{1}{2}} = Q^n_{j+1} - \frac{\Delta x}{2}(Q_x)^n_{j+1}.
\end{equation*}
which are considered to define a Riemann problem at $x= x_\jph$. The solution of the Riemann problems gives rise to some waves possibly leading to non-smooth solution. We separate the spatial domain into a smooth part and non-smooth part by estimating the maximal local wave speeds $a^{\pm}$ at interfaces;
\begin{equation*}
    a^-_{j+\frac{1}{2}} = \min \{\lambda_1(Q^-_{j+\frac{1}{2}}), \lambda_1(Q^+_{j+\frac{1}{2}}),0 \}\leq 0, \quad
    a^+_{j+\frac{1}{2}} = \max \{\lambda_d(Q^-_{j+\frac{1}{2}}), \lambda_d(Q^+_{j+\frac{1}{2}}),0 \}\geq 0,
\end{equation*}
where $\lambda$ is the eigenvalue of $\frac{\partial f(Q)}{\partial Q}$ and $\lambda_1\leq\lambda_2\leq...\leq\lambda_d$. With reference to Figure \ref{unORsmooth}, the solution is smooth in the green part which is the interval $[x_{j-
\frac{1}{2},r}, x_{j+\frac{1}{2},l}]$, with  
$x_{j+\frac{1}{2},l} :=x_\jph + a^-_\jph \Delta t$ and $ x_{\jph,r} := x_\jph + a^+_\jph \Delta t$, and possibly non-smooth in the red portions where the waves arising from the Riemann problems are present.

To obtain the intermediate values $Q^{\text{int}}_j$ and $Q^{\text{int}}_\jph$, we integrate the system \eqref{1Dlaw} over the green part $[x_{\jmh,r}, x_{\jph,l}]\times[t^n, t^{n+1}]$ and the red part $[x_{\jph,l}, x_{\jph,r}]\times[t^n, t^{n+1}]$, respectively. 
\begin{align*}
    Q^{\text{int}}_j 
    :=&\frac{1}{\Delta x^n_j}\int^{x_{\jph,l}}_{x_{\jmh,r}} Q(x,t^{n+1})dx \\
    =& \frac{1}{\Delta x^n_j}\Bigg[\int^{x_{\jph,l}}_{x_{\jmh,r}} Q(x,t^n)dx - \int^{t^{n+1}}_{t}\Big(f(Q(x_{\jph,l},t))-f(Q(x_{\jmh,r},t))\Big)\,dt\Bigg],\\
    Q^{\text{int}}_\jph 
    :=&\frac{1}{\Delta x^n_\jph}\int^{x_{\jph,r}}_{x_{\jph,l}} Q(x,t^{n+1})dx \\
    =& \frac{1}{\Delta x^n_\jph}\Bigg[\int^{x_{\jph,r}}_{x_{\jph,l}} Q(x,t^n)dx - \int^{t^{n+1}}_{t}\Big(f(Q(x_{\jph,r},t))-f(Q(x_{\jph,l},t))\Big)\,dt\Bigg],  
\end{align*}
where $\Delta x^n_j=x_{\jph,l}-x_{\jmh,r}$ and $\Delta x^n_\jph=x_{\jph,r}-x_{\jph,l}$ are the length of the interval.

Applying reconstruction \eqref{1DRec} to the first integral and the midpoint rule to the second integral, we rewrite the above approximation 
\begin{align*}
    Q^{\text{int}}_j &= Q^n_j+(x_{j,m}-x_j)(Q_x)^n_j - \frac{\Delta t}{\Delta x^n_j}\bigg[f(Q^{n+\frac{1}{2}}_{\jph,l})-f(Q^{n+\frac{1}{2}}_{\jmh,r})\bigg],  \\
    Q^{\text{int}}_\jph &= 
    \frac{1}{\Delta x^n_\jph}\times\\
    &\times\bigg[-a^n_\jph\Delta t \Big(Q^n_j+(x_{\jph,lm}-x_j)(Q_x)^n_j \Big) +a^+_\jph\Delta t \Big(Q^n_{j+1}+(x_{\jph,rm}-x_{j+1})(Q_x)^n_{j+1} \Big) \bigg] \\
    & - \frac{\Delta t}{\Delta x^n_\jph}\bigg[f(Q^{n+\frac{1}{2}}_{\jph,r})-f(Q^{n+\frac{1}{2}}_{\jph,l})\bigg].
\end{align*}
where $x_{j,m}=(x_{\jph,l}+x_{\jmh,r})/2$ and $x_{\jph,lm(rm)}=(x_{\jph,l(r)}+x_\jph)/2$ denote the midpoint of the corresponding interval. The second order approximations of $Q^{n+\frac{1}{2}}_{\jph,r(l)}$ are obtained by Taylor expansion,
\begin{align*}
    Q^{n+\half}_{\jph,l}&:=Q^n_{\jph,l}+\frac{\Delta t}{2}(Q_t)^n_{\jph,l} \approx Q^n_{\jph,l} - \frac{\Delta t}{2}[f(Q)_x]^n_j, \\
    Q^{n+\half}_{\jph,r}&:=Q^n_{\jph,r}+\frac{\Delta t}{2}(Q_t)^n_{\jph,r} \approx Q^n_{\jph,r} - \frac{\Delta t}{2}[f(Q)_x]^n_{j+1},
\end{align*}
with
\begin{align*}
    Q^n_{\jph,l}:=Q^n_j+(\frac{\Delta x}{2}+a^-_\jph \Delta t)(Q_x)^n_j,\quad
    Q^n_{\jph,r}:=Q^n_{j+1}-(\frac{\Delta x}{2}-a^+_\jph \Delta t)(Q_x)^n_{j+1}. 
\end{align*}
~\\
\noindent\textbf{Step 3: Projection:}
In the last step, we project back to the cell averages on each  cell $C_j$. Differently from the previous central scheme using a linear approximation to reconstruct $Q^{\text{int}}_\jph$ over Riemann fans, the piecewise constant is considered here.  

In order to represent the jump in density across the contact wave, we split the intermediate values $Q^{\text{int}}_\jph$ into two constant values at two sides of the contact, which can make the largest difference between the two sides, but without creating new local extrema. 
\begin{align*}
    Q^\text{int}_\jph(x) = \left\{
    \begin{array}{cc}
         & Q^L_\jph,\qquad x_{\jph,l}<x<x_\jph+v^\text{int}_\jph\Delta t,   \\
         & Q^R_\jph,\qquad x_\jph+v^\text{int}_\jph\Delta t<x<x_{\jph,r}. 
    \end{array}
    \right.
\end{align*}

Following the algorithm in~\cite{LDCU2} and considering the velocity along $x$ axis, $v_1$, we reconstruct the density $\rho$ by
\begin{equation} \label{1DMHD_rhoRec}
    \rho^L_\jph = \rho^\text{int}_\jph + \frac{\delta_\jph}{a^-_\jph+(v_1)^\text{int}_\jph},\quad\text{and}\quad 
    \rho^R_\jph = \rho^\text{int}_\jph + \frac{\delta_\jph}{a^+_\jph+(v_1)^\text{int}_\jph},
\end{equation}
with 
\begin{equation*} 
    \delta_\jph:= \text{minmod} (S^-_\jph, S^+_\jph),
\end{equation*}
where $S^\pm$ are defined by
\begin{equation*} 
    S^-_\jph := \Big((v_1)^\text{int}_\jph-a^-_\jph\Big) \Big(\rho^\text{int}_\jph-\rho^\text{int}_{\jph,l}\Big), \quad 
    S^+_\jph := \Big(a^+_\jph-(v_1)^\text{int}_\jph\Big) \Big(\rho^\text{int}_{\jph,l}-\rho^\text{int}_\jph\Big).
\end{equation*}
Next, due to the properties that velocity $v_i$ and magnetic field $B_i$ are continuous at the contact discontinuity while $B_1\neq0$, we define momentum and magnetic field at the right and left-hand sides of the contact as follows,
\begin{align} 
    & \text{Momentum:} \qquad (\rho v_i)^{L(R)}_\jph=\rho^{L(R)}_\jph (v_i)^\text{int}_\jph, \qquad  i=\{1,2,3\},  \label{1DMHD_moRec}\\
    & \text{Magnetic field:} \;\qquad B_i^{L(R)}=B^\text{int}_i, \qquad \qquad \qquad \,\;\; i=\{2,3\}. \label{1DMHD_BRec}
\end{align}
Here, $(v_i)^\text{int}_\jph=\frac{(\rho v_i)^\text{int}_\jph}{\rho^\text{int}_\jph}$.  
Since the pressure $p$ is also continuous at contact discontinuity and according to the definition of energy, $E=\frac{p}{\gamma-1}+ \frac{1}{2} \rho|\textbf{V}|^2 +\frac{1}{2}|\textbf{B}|^2$, we obtain 
\begin{align} 
    E^\text{int}_\jph=&\frac{p^\text{int}_\jph}{\gamma-1}+ \frac{1}{2} \rho^\text{int}_\jph|\textbf{V}^\text{int}_\jph|^2 +\frac{1}{2}|\textbf{B}^\text{int}_\jph|^2, \label{1DE_int}    \\
    E^{L(R)}_\jph=&\frac{p^\text{int}_\jph}{\gamma-1}+ \frac{1}{2} \rho^{L(R)}_\jph|\textbf{V}^\text{int}_\jph|^2 +\frac{1}{2}|\textbf{B}^\text{int}_\jph|^2 \label{1DMHD_E_LR}.
\end{align}
Then, substituting Eq.~\eqref{1DE_int} into Eq.~\eqref{1DMHD_E_LR}, the reconstructed energy can be written as
\begin{align} \label{1DMHD_ERec}
    E^{L(R)}_\jph = E^\text{int}_\jph +\frac{\rho^{L(R)}_\jph - \rho^\text{int}_\jph}{2}|\textbf{V}^\text{int}_\jph|^2, 
\end{align}
Next, to project back to the original unstaggered cell $[x_{j-\frac{1}{2}}, x_{j+\frac{1}{2}}]$, we substitute the above reconstructed values \eqref{1DMHD_rhoRec}, \eqref{1DMHD_moRec}, \eqref{1DMHD_BRec}, and \eqref{1DMHD_ERec} into 
\begin{align}
\label{ldcu1dup}
    Q^{n+1}_j := \frac{1}{\Delta x}\int^{x_\jph}_{x_\jmh} Q^{\text{int}}(x)dx,    
\end{align}
and finally obtain the quantities in the next time step $t^{n+1}$ with the formula 
\begin{equation} \label{1DMHDfully1}
\begin{split}
    Q^{n+1} = Q^\text{int}_j + 
    \frac{\Delta t}{\Delta x}\bigg[ &a^+_\jmh(Q^R_\jmh-Q^\text{int}_j) +\max((v_1)^\text{int}_\jmh,0)\big( Q^L_\jmh-Q^R_\jmh \big) \\ 
    &-a^-_\jph(Q^L_\jph-Q^\text{int}_j)+\min((v_1)^\text{int}_\jph,0)\big( Q^L_\jph - Q^R_\jph \big) \bigg]
\end{split}
\end{equation}

\begin{remark} \label{remark}
Since the reconstruction in the projection step is piecewise constant, the  resulting scheme is first-order accurate.  \end{remark}

\subsection{1-D LDCU MHD semi-discrete scheme} \label{1DMHD}
In this section, we construct a semi-discrete scheme from the 1-D LDCU MHD fully-discrete scheme in previous section.  Firstly, we rewrite the formula \eqref{1DMHDfully1} in a form as
\begin{align} \label{1DMHDfully2}
    Q^{n+1}_j 
    =Q^\text{int}_j + \frac{\Delta t}{\Delta x} \bigg[a^+_\jmh(Q^\text{int}_\jmh-Q^\text{int}_j) - a^-_\jph(Q^\text{int}_\jph-Q^\text{int}_j) + \alpha_\jmh\hat{\delta}_\jmh - \alpha_\jph\hat{\delta}_\jph \bigg].
\end{align} 
where, $\alpha_\jph\hat{\delta}_\jph=:\mathrm{K}^\text{int}_\jph$ is the correction term produced by the LDCU method to reduce dissipation, and is defined by 
\begin{align*}
    \alpha_{j+\frac{1}{2}}=
    \begin{cases}
         \dfrac{a^-_\jph}{a^-_\jph-(v_1)^\text{int}_\jph},\qquad (v_1)^\text{int}_\jph\geq 0,  &  \\
        \dfrac{a^+_\jph}{a^+_\jph-(v_1)^\text{int}_\jph},\qquad (v_1)^\text{int}_\jph< 0, &
    \end{cases}
    \qquad
    \hat{\delta}_{j+\frac{1}{2}} =
    \delta_{j+\frac{1}{2}}
    \begin{bmatrix}
          1 \\
          (v_1)^\text{int}_{j+\frac{1}{2}} \\
          (v_2)^\text{int}_{j+\frac{1}{2}}\\
          (v_3)^\text{int}_{j+\frac{1}{2}}\\
          0 \\
          0 \\
          \frac{1}{2} |\textbf{V}^\text{int}_{j+\frac{1}{2}}|^2
    \end{bmatrix},
\end{align*}
with
\begin{align*}
    \delta_{j+\frac{1}{2}}= \text{minmod}\bigg(\Big((v_1)^\text{int}_{j+\frac{1}{2}}-a^-_{j+\frac{1}{2}}\Big) \Big(\rho^\text{int}_{j+\frac{1}{2}}-\rho^\text{int}_{j+\frac{1}{2},l}\Big),\; \Big(a^+_{j+\frac{1}{2}}-(v_1)^\text{int}_{j+\frac{1}{2}}\Big)\Big(\rho^\text{int}_{j+\frac{1}{2},r}-\rho^\text{int}_{j+\frac{1}{2}}\Big)\bigg). 
\end{align*}
Then, we use $Q^{n+1}_j$ in formula \eqref{1DMHDfully2} to compute the difference $Q^{n+1}_j-Q^n_j$ and divide it by the time deviation $\Delta t$, and next, take the limit of it, i.e.,
\begin{equation*}
    \lim_{\Delta t \to 0} \frac{Q^{n+1}_j-Q^n_j}{\Delta t}.
\end{equation*}
After some calculations, we obtain the semi-discrete scheme
\begin{align} \label{1Dsemi}
    \frac{dQ_j(t)}{dt} = -\frac{F_{j+\frac{1}{2}}(t)-F_{j-\frac{1}{2}}(t)}{\Delta x},
\end{align}
with the numerical flux, which is defined by
\begin{align} \label{1Dflux}
    F_{j+\frac{1}{2}}(t) 
    =\frac{a^+_{j+\frac{1}{2}}a^-_{j-\frac{1}{2}}}{a^+_{j+\frac{1}{2}}-a^-_{j+\frac{1}{2}}}\Big(Q^+_{j+\frac{1}{2}}-Q^-_{j+\frac{1}{2}}\Big) + \frac{\bigg[a^+_{j+\frac{1}{2}}f(Q^-_{j+\frac{1}{2}})-a^-_{j+\frac{1}{2}}f(Q^+_{j+\frac{1}{2}})\bigg]}{a^+_{j+\frac{1}{2}}-a^-_{j+\frac{1}{2}}} + \mathrm{K}^*_{j+\frac{1}{2}},
\end{align}
and the correction term $\mathrm{K}^*$ is given by
\begin{align*}
    \mathrm{K}^*_\jph = \alpha^*_\jph\delta^*_\jph
    \begin{bmatrix}
        1 \\
        (v_1)^*_\jph \\
        (v_2)^*_\jph\\
        (v_3)^*_\jph\\
        0 \\
        0 \\
        \frac{1}{2} |{\mathbf{V}}^*_\jph|^2
    \end{bmatrix},
    \quad \text{where} \quad
    \alpha^*_\jph=
    \begin{cases}
         \dfrac{a^-_\jph}{a^-_\jph-(v_1)^*_\jph}, & (v_1)^*_\jph \geq 0, \\
         \dfrac{a^+_\jph}{a^+_\jph-(v_1)^*_\jph}, & (v_1)^*_\jph < 0,
    \end{cases}
\end{align*}
and
\begin{align*}
    \delta^*_{j+\frac{q}{2}} = \text{minmod}\bigg(\Big((v_1)^*_{j+\frac{1}{2}}-a^-_{j+\frac{1}{2}}\Big) \Big(\rho^*_{j+\frac{1}{2}}-\rho^-_{j+\frac{1}{2}}\Big), \; \Big(a^+_{j+\frac{1}{2}}-(v_1)^*_{j+\frac{1}{2}}\Big) \Big(\rho^+_{j+\frac{1}{2}}-\rho^*_{j+\frac{1}{2}}\Big)\bigg),
\end{align*}
and $\rho^*_{j+\frac{1}{2}}$ and $(v_i)^*_{j+\frac{1}{2}}:=\frac{(\rho v_i)^*_{j+\frac{1}{2}}}{\rho^*_{j+\frac{1}{2}}}$ are obtained by
\begin{align}\label{q_star}
    Q^*_{j+\frac{1}{2}}
    :=&\lim_{\Delta t \to 0} Q^\text{int}_{j+\frac{1}{2}.} \nonumber \\
    =&\frac{a^+_{j+\frac{1}{2}}Q^+_{j+\frac{1}{2}}-a^-_{j+\frac{1}{2}}Q^-_{j+\frac{1}{2}}-\Big[f(Q^+_{j+\frac{1}{2}})-f(Q^-_{j+\frac{1}{2}})\Big]}{a^+_{j+\frac{1}{2}}-a^-_{j+\frac{1}{2}}}.
\end{align}

\section{2-D LDCU MHD scheme}
We now extend the  1-D semi-discrete scheme from the previous section to the 2-D system by using the same idea of improving the resolution of contact discontinuity. Consider the 2-D MHD equations,
\begin{equation*}
    q_t+f(q)_x+g(q)_y =0,
\end{equation*}
and
\begin{equation*}
    q = \left[
    \begin{matrix}
        \rho  \\
        \rho v_1 \\
        \rho v_2 \\
        \rho v_3 \\
        B_1 \\
        B_2 \\
        B_3 \\
        E
    \end{matrix}
    \right],\quad
    f(q) = \left[
    \begin{matrix}
        \rho v_1  \\
        \rho v_1^2 + p + \frac{1}{2}|{\bf{B}}|^2 - B_1^2 \\
        \rho v_1v_2 - B_1B_2 \\
        \rho v_1v_3 - B_1B_3 \\
        0 \\
        B_2v_1 - B1v_2 \\
        B_3v_1 - B_1v_3 \\
        (E+p+\frac{1}{2}|{\bf{B}}|^2)v_1 - B_1(B_1v_1+B_2v_2+B_3v_3)
    \end{matrix}
    \right],
\end{equation*}
\begin{equation*}
    g(q) = \left[
    \begin{matrix}
        \rho v_2  \\
        \rho v_2v_1 - B_2B_1 \\
        \rho v_2^2 + p + \frac{1}{2}|{\bf{B}}|^2 - B_2^2 \\
        \rho v_2v_3 - B_2B_3 \\
        B_1v_2 - B_2v_1 \\
        0 \\
        B_3v_2 - B_2v_3 \\
        (E+p+\frac{1}{2}|{\bf{B}}|^2)v_2 - B_2(B_1v_1+B_2v_2+B_3v_3)
    \end{matrix}
    \right].
\end{equation*}
We separate the eight variables into two groups:
\begin{itemize}
\item Hydrodynamic variables $U=(\rho, \rho v_1, \rho v_2, \rho v_3, E)^\top$
\item  Magnetic variables $\mathbf{B} = (B_1, B_2, B_3)^\top$
\end{itemize}
The hydrodynamics variables are stored at the cell centers while the magnetic variables $(B_1,  B_2)$ are stored on the cell faces. For the hydrodynamic variables $U$, we use the 'dimension-by-dimension' method to extend the 1-D semi-discrete LDCU MHD scheme in Section \ref{1DMHD} to 2-D, and for the magnetic variables \textbf{B}, the upwind constrained transport HLL method is used with a time evolution scheme, which will be discussed in Section \ref{UCTHLL}. 

Before we move to the next section, we need to define two things first: a second-order linear reconstruction with the cell-averaged 
$Q_{j,k}$  
\begin{align}
    Q^+_{j+\frac{1}{2},k} &:= Q_{j+1,k} - \frac{\Delta x}{2}(Q_x)_{j+1,k},\qquad 
    Q^+_{j,k+\frac{1}{2}} := Q_{j,k+1} - \frac{\Delta y}{2}(Q_y)_{j,k+1}, \label{RQ+} \\
    Q^-_{j+\frac{1}{2},k} &:= Q_{j,k} + \frac{\Delta x}{2}(Q_x)_{j,k},\qquad\qquad\,
    Q^-_{j,k+\frac{1}{2}} := Q_{j,k} + \frac{\Delta y}{2}(Q_y)_{j,k}, \label{RQ-}
\end{align}
and the maximal wave speeds $a^\pm$ and $b^\pm$ at interfaces 
\begin{align*}
    a^-_{j+\frac{1}{2},k} =& \min \{\lambda_1(Q^-_{j+\frac{1}{2},k}), \lambda_1(Q^+_{j+\frac{1}{2},k}),-10^{-8} \}, \quad
    a^+_{j+\frac{1}{2},k} = \max \{\lambda_d(Q^-_{j+\frac{1}{2},k}), \lambda_d(Q^+_{j+\frac{1}{2},k}),10^{-8} \}, \\
    b^-_{j,k+\frac{1}{2}} =& \min \{\lambda_1(Q^-_{j,k+\frac{1}{2}}), \lambda_1(Q^+_{j,k+\frac{1}{2}}),-10^{-8} \}, \quad\:
    b^+_{j,k+\frac{1}{2}} = \max \{\lambda_d(Q^-_{j,k+\frac{1}{2}}), \lambda_d(Q^+_{j,k+\frac{1}{2}}),10^{-8} \},
\end{align*}
where $\lambda$ are the eigenvalues of the flux Jacobian.

\subsection{2-D LDCU scheme}
The domain is assumed to be rectangular and is partitioned into disjoint cells $\Cjk=(x_\jmh, x_\jph) \times (y_\kmh, y_\kph)$ of width $(\Delta x, \Delta y)$ along the two axes. Let $(x_j,y_k)$ denote the center of cell $\Cjk$ and $Q_{j,k}$ denote the cell average solution. 

Consider the hydrodynamic variables $U$. We extend the 1-D semi-discrete MHD scheme \eqref{1Dsemi} to 2-D to update only the $U$ variables.  Therefore, the 2-D semi-discrete scheme is in the form of 
\begin{align} \label{semiQ}
    \frac{dU_{j,k}(t)}{dt} = -\frac{F_{j+\frac{1}{2},k}(t)-F_{j-\frac{1}{2},k}(t)}{\Delta x} - \frac{G_{j,k+\frac{1}{2}}(t)-G_{j,k-\frac{1}{2}}(t)}{\Delta y},
\end{align}
with the numerical fluxes
\begin{align*}
    F_{j+\frac{1}{2},k}(t) 
    =\frac{a^+_{j+\frac{1}{2},k}a^-_{j+\frac{1}{2},k}}{a^+_{j+\frac{1}{2},k}-a^-_{j+\frac{1}{2},k}}\Big(U^+_{j+\frac{1}{2},k}-U^-_{j+\frac{1}{2},k}\Big) + \frac{\bigg[a^+_{j+\frac{1}{2},k}f(U^-_{j+\frac{1}{2},k})-a^-_{j+\frac{1}{2},k}f(U^+_{j+\frac{1}{2},k})\bigg]}{a^+_{j+\frac{1}{2},k}-a^-_{j+\frac{1}{2},k}} + {\mathrm{K}}^*_{j+\frac{1}{2},k},
\end{align*}
\begin{align*}
    G_{j,k+\frac{1}{2}}(t) 
    =\frac{b^+_{j,k+\frac{1}{2}}b^-_{j,k+\frac{1}{2}}}{b^+_{j,k+\frac{1}{2}}-b^-_{j,k+\frac{1}{2}}}\Big(U^+_{j,k+\frac{1}{2}}-U^-_{j,k+\frac{1}{2}}\Big) + \frac{\bigg[b^+_{j,k+\frac{1}{2}}g(U^-_{j,k+\frac{1}{2}})-b^-_{j,k+\frac{1}{2}}g(U^+_{j,k+\frac{1}{2}})\bigg]}{b^+_{j,k+\frac{1}{2}}-b^-_{j,k+\frac{1}{2}}} + {\mathrm{K}}^*_{j,k+\frac{1}{2}},
\end{align*}
where $U^\pm$ are defined in Eq~\eqref{RQ+} and \eqref{RQ-}.  Similar to the 1-D scheme, the correction terms from LDCU scheme, $\mathrm{K}^*$,
are defined by
\begin{align*}
    \mathrm{K}^*_{j+\frac{1}{2},k} = \alpha^*_{j+\frac{1}{2},k}\delta^*_{j+\frac{1}{2},k}
    \left[
    \begin{matrix}
        1 \\
        (v_1)^*_{j+\frac{1}{2},k} \\
        (v_2)^*_{j+\frac{1}{2},k}\\
        (v_3)^*_{j+\frac{1}{2},k}\\
        \frac{1}{2} |\textbf{V}^*_{j+\frac{1}{2},k}|^2
    \end{matrix}
    \right],
\end{align*}
where
\begin{equation*}
\alpha^*_{j+\frac{1}{2},k}=
\left \{
\begin{array}{cc}
     \dfrac{a^-_{j+\frac{1}{2},k}}{a^-_{j+\frac{1}{2},k}-(v_1)^*_{j+\frac{1}{2},k}},\qquad (v_1)^*_{j+\frac{1}{2},k}\geq 0,  &  \\
     \dfrac{a^+_{j+\frac{1}{2},k}}{a^+_{j+\frac{1}{2},k}-(v_1)^*_{j+\frac{1}{2},k}},\qquad (v_1)^*_{j+\frac{1}{2},k}< 0, & 
\end{array}
\right.
\end{equation*}
\begin{align*}
    \delta^*_{j+\frac{q}{2},k} = \text{minmod}\bigg(\Big((v_1)^*_{j+\frac{1}{2},k}-a^-_{j+\frac{1}{2},k}\Big) \Big(\rho^*_{j+\frac{1}{2},k}-\rho^-_{j+\frac{1}{2},k}\Big), \Big(a^+_{j+\frac{1}{2},k}-(v_1)^*_{j+\frac{1}{2},k}\Big) \Big(\rho^+_{j+\frac{1}{2},k}-\rho^*_{j+\frac{1}{2},k}\Big)\bigg),
\end{align*}
and $\rho^*$ and $(v_i)^*$ are obtained by  
\begin{align} \label{q_star_x}
    Q^*_{j+\frac{1}{2},k}
    =\frac{a^+_{j+\frac{1}{2},k}Q^+_{j+\frac{1}{2},k}-a^-_{j+\frac{1}{2},k}Q^-_{j+\frac{1}{2},k}-\left[f(Q^+_{j+\frac{1}{2},k})-f(Q^-_{j+\frac{1}{2},k})\right]}{a^+_{j+\frac{1}{2},k}-a^-_{j+\frac{1}{2},k}},
\end{align}
which are similar to Eq.~\eqref{q_star}, while $\mathrm{K}^*_{j,k+\frac{1}{2}}$ is defined the same way that 
\begin{align*}
    \mathrm{K}^*_{j,k+\frac{1}{2}} = \alpha^*_{j,k+\frac{1}{2}}\delta^*_{j,k+\frac{1}{2}}
    \left[
    \begin{matrix}
        1 \\
        (v_1)^*_{j,k+\frac{1}{2}} \\
        (v_2)^*_{j,k+\frac{1}{2}}\\
        (v_3)^*_{j,k+\frac{1}{2}}\\
        \frac{1}{2} |\textbf{V}^*_{j,k+\frac{1}{2}}|^2
    \end{matrix}
    \right],
\end{align*}
where
\begin{equation*}
\alpha^*_{j,k+\frac{1}{2}}=
\left \{
\begin{array}{cc}
     \dfrac{b^-_{j,k+\frac{1}{2}}}{b^-_{j,k+\frac{1}{2}}-(v_2)^*_{j,k+\frac{1}{2}}},\qquad (v_2)^*_{j,k+\frac{1}{2}}\geq 0,  &  \\
     \dfrac{b^+_{j,k+\frac{1}{2}}}{b^+_{j,k+\frac{1}{2}}-(v_2)^*_{j,k+\frac{1}{2}}},\qquad (v_2)^*_{j,k+\frac{1}{2}}< 0, & 
\end{array}
\right.
\end{equation*}
\begin{align*}
    \delta^*_{j,k+\frac{1}{2}} = \text{minmod}\bigg(\Big((v_2)^*_{j,k+\frac{1}{2}}-b^-_{j,k+\frac{1}{2}}\Big) \Big(\rho^*_{j,k+\frac{1}{2}}-\rho^-_{j,k+\frac{1}{2}}\Big),\; \Big(b^+_{j,k+\frac{1}{2}}-(v_2)^*_{j,k+\frac{1}{2}}\Big)\Big(\rho^+_{j,k+\frac{1}{2}}-\rho^*_{j,k+\frac{1}{2}}\Big)\bigg),
\end{align*}   
and
\begin{align} \label{q_star_y}
    Q^*_{j,k+\frac{1}{2}}
    =&\frac{b^+_{j,k+\frac{1}{2}}Q^+_{j,k+\frac{1}{2}}-b^-_{j,k+\frac{1}{2}}Q^-_{j,k+\frac{1}{2}}-\Big[g(Q^+_{j,k+\frac{1}{2}})-g(Q^-_{j,k+\frac{1}{2}})\Big]}{b^+_{j,k+\frac{1}{2}}-b^-_{j,k+\frac{1}{2}}}.
\end{align}
\\
\begin{remark}
Since we only discuss the 2-D experiments in this paper, the magnetic field component along $z$-axis denoted by $B_3$ does not need to be involved in the divergence-free method. It can simply use the LDCU scheme in the evolution step, that is we consider $U=(\rho, \rho v_1, \rho v_2, \rho v_3, B_3, E)^\top$ with $\mathrm{K}^*[B_3]=0$, which is similar to the 1-D LDCU MHD scheme.
\end{remark}

\subsection{Upwind constrained transport HLL method}
\label{UCTHLL}
To ensure that the divergence of magnetic field is zero, we use the UCT-HLL method from Eq.~(28) and Eq.~(29) in Section 4.2 of \cite{CT}. Consider the semi-discrete form of the induction equation
\begin{align}
    \frac{d(B_1)_{j+\frac{1}{2},k}}{dt} 
    &= -\frac{\Omega_{j+\frac{1}{2},k+\frac{1}{2}}-\Omega_{j+\frac{1}{2},k-\frac{1}{2}}}{\Delta y}  \label{semiB1}, \\
    \frac{d(B_2)_{j,k+\frac{1}{2}}}{dt} 
    &=+\frac{\Omega_{j+\frac{1}{2},k+\frac{1}{2}}+\Omega_{j-\frac{1}{2},k+\frac{1}{2}}}{\Delta x}  \label{semiB2},
\end{align}
where the electric field $\Omega$ at the nodes $(\jph,\kph)$ is defined by
\begin{align} \label{EF}
\begin{aligned}
    \Omega_{j+\frac{1}{2},k+\frac{1}{2}} 
    = - &\frac{\alpha^+_x(V_1 B_2)^W_{j+\frac{1}{2},k+\frac{1}{2}}+\alpha^-_x(V_1 B_2)^E_{j+\frac{1}{2},k+\frac{1}{2}} - \alpha^+_x\alpha^-_x((B_2)^E_{j+\frac{1}{2},k+\frac{1}{2}}-(B_2)^W_{j+\frac{1}{2},k+\frac{1}{2}})}{\alpha^+_x + \alpha^-_x} \\
    + &\frac{\alpha^+_y(V_2 B_1)^S_{j+\frac{1}{2},k+\frac{1}{2}}+\alpha^-_y(V_2 B_1)^N_{j+\frac{1}{2},k+\frac{1}{2}} - \alpha^+_y\alpha^-_y((B_1)^N_{j+\frac{1}{2},k+\frac{1}{2}}-(B_1)^S_{j+\frac{1}{2},k+\frac{1}{2}})}{\alpha^+_y + \alpha^-_y},
\end{aligned}
\end{align}
\\
Here, $\alpha^{\pm}$ are defined by
\begin{align*}
\begin{aligned}
    \alpha^+_x = \max(0, \lambda^R_{j+\frac{1}{2},k}, \lambda^R_{j+\frac{1}{2},k+1}), \quad
    \alpha^-_x = -\min(0, \lambda^L_{j+\frac{1}{2},k}, \lambda^L_{j+\frac{1}{2},k+1}), \\
    \alpha^+_y = \max(0, \lambda^R_{j,k+\frac{1}{2}}, \lambda^R_{j+1,k+\frac{1}{2}}), \quad
    \alpha^-_y = -\min(0, \lambda^L_{j,k+\frac{1}{2}}, \lambda^L_{j+1,k+\frac{1}{2}}).
\end{aligned}
\end{align*}
For transverse velocities $V$, we first compute the velocities at the interface, e.g., $\overline{V}_2$ at $x$-interface,  
\begin{align*}
    (\overline V_2)_{j+\frac{1}{2},k} = \frac{\alpha^R_x
     (v_2)^-_{j+\frac{1}{2},k} + \alpha^L_x (v_2)^+_{j+\frac{1}{2}.k}}{\alpha^R_x +\alpha^L_x},
\end{align*}
with
\begin{align*}
    \alpha^R_x = \max(0, \lambda^R_{j+\frac{1}{2},k}),\quad
    \alpha^L_x = -\min(0,\lambda^L_{j+\frac{1}{2},k}),
\end{align*}
and $(v_2)^\pm_{j+\frac{1}{2},k}$ are obtained from Eq.~\eqref{RQ+} and \eqref{RQ-}, and then reconstruct the velocities from the face to the corner along the transverse direction with a proper limiter $\phi$, (refer to Appendix B of~\cite{DefV})
\begin{align} \label{RV}
    (V_1)^W_{j+\frac{1}{2},k+\frac{1}{2}} = (\overline{V}_1)_{j,k+\frac{1}{2}} + \frac{\Delta x}{2}(\phi_x)_{j,k+\frac{1}{2}},\quad
    (V_1)^E_{j+\frac{1}{2},k+\frac{1}{2}} = (\overline{V}_1)_{j+1,k+\frac{1}{2}} - \frac{\Delta x}{2}(\phi_x)_{j+1,k+\frac{1}{2}},\\
    (V_2)^S_{j+\frac{1}{2},k+\frac{1}{2}} = (\overline{V}_2)_{j+\frac{1}{2},k} + \frac{\Delta y}{2}(\phi_y)_{j+\frac{1}{2},k},\quad
    (V_2)^N_{j+\frac{1}{2},k+\frac{1}{2}} = (\overline{V}_2)_{j+\frac{1}{2},k+1} - \frac{\Delta y}{2}(\phi_y)_{j+\frac{1}{2},k+1},
\end{align}
where
\begin{align} 
    (\phi_x)_{j,k+\frac{1}{2}} = \text{(MC-$\theta$)-limiter}((\overline{V}_1)_{j-1,k+\frac{1}{2}},\; (\overline{V}_1)_{j,k+\frac{1}{2}},\;(\overline{V}_1)_{j+1,k+\frac{1}{2}}),\\
    (\phi_y)_{j+\frac{1}{2},k} = \text{(MC-$\theta$)-limiter}((\overline{V}_2)_{j+\frac{1}{2},k-1},\; (\overline{V}_2)_{j+\frac{1}{2},k},\;(\overline{V}_2)_{j+\frac{1}{2},k+1}). \label{dV}
\end{align}
Similarly, we can use the same method as in Eq.~\eqref{RV}-\eqref{dV}, to reconstruct magnetic fields $B_2$ and $B_1$ from the interface to the corner.  Then, we apply the Runge-kutta method to obtain the magnetic field $(B_1)_{j+\frac{1}{2},k}$ and $(B_2)_{j,k+\frac{1}{2}}$ at the next time step. Finally, the magnetic field at the cell center can be obtained by  
\begin{align*}
    (B_1)^{n+1}_{j,k} = \frac{(B_1)^{n+1}_{j+\frac{1}{2},k}+(B_1)^{n+1}_{j-\frac{1}{2},k}}{2}
    \quad\text{and}\quad
    (B_2)^{n+1}_{j,k} = \frac{(B_2)^{n+1}_{j,k+\frac{1}{2}}+(B_2)^{n+1}_{j,k-\frac{1}{2}}}{2}.
\end{align*}
By construction, this method maintains the constancy of the magnetic field divergence measured at the cell centers, and this is shown in the 2-D test cases.

\section{Numerical results}\label{Numerical}
To demonstrate our new scheme, we apply our 1-D fully-discrete scheme \eqref{1DMHDfully1} and 2-D semi-discrete scheme \eqref{semiQ}, \eqref{semiB1}, \eqref{semiB2} to several numerical experiments.  Furthermore, we  compare the results with the scheme without the LDCU correction to show the difference.

\subsection{1-D experiments}
In the 1-D tests, we use the following formula to compute the time step $\Delta t$, and the CFL number is set as 0.4.   
\begin{align*}
    \Delta t = \text{CFL} \cdot \min_{j} \frac{\Delta x}{\max(a^+_{j+\frac{1}{2}}, a^-_{j+\frac{1}{2}})}.
\end{align*}
The boundary condition used for these 1-D experiments is the outflow boundary condition. The solution we compared with is using the scheme without the reconstruction $Q^{L(R)}_\jph$ in the projection step, which means the intermediate average $Q^\text{int}_\jph$ is used to project back to the uniform cell.  The exact solutions are obtained by the compared scheme using a finer grid with 6000 points. The obtained 1-D results are first-order accurate due to the reason in Remark~\ref{remark}.

\subsubsection{Brio-Wu shock tube problem} \label{sec.1DBW}
\begin{figure}[htbp]
    \centering 
    \subfigure[]{\includegraphics[height=0.36\linewidth]{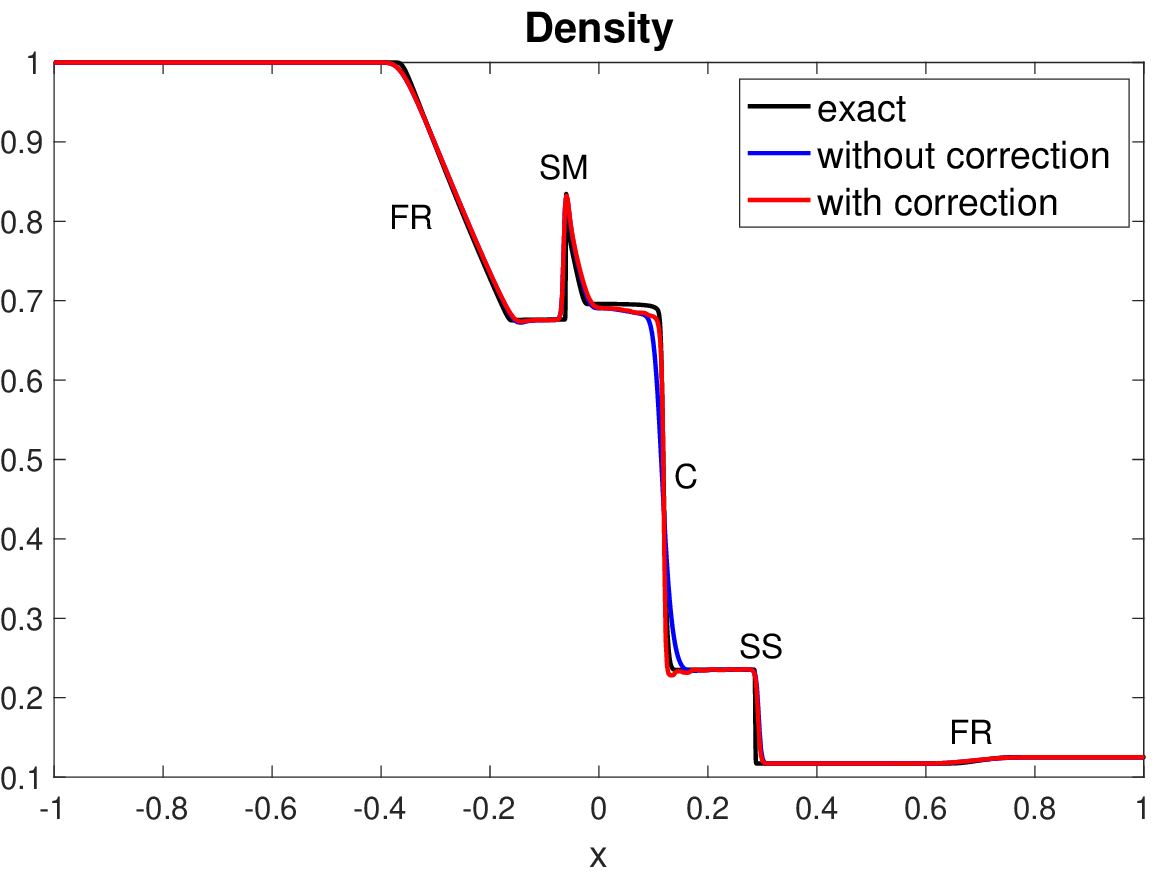}
    }
    \subfigure[Zoomed view around contact wave]{
    \includegraphics[height=0.36\linewidth]{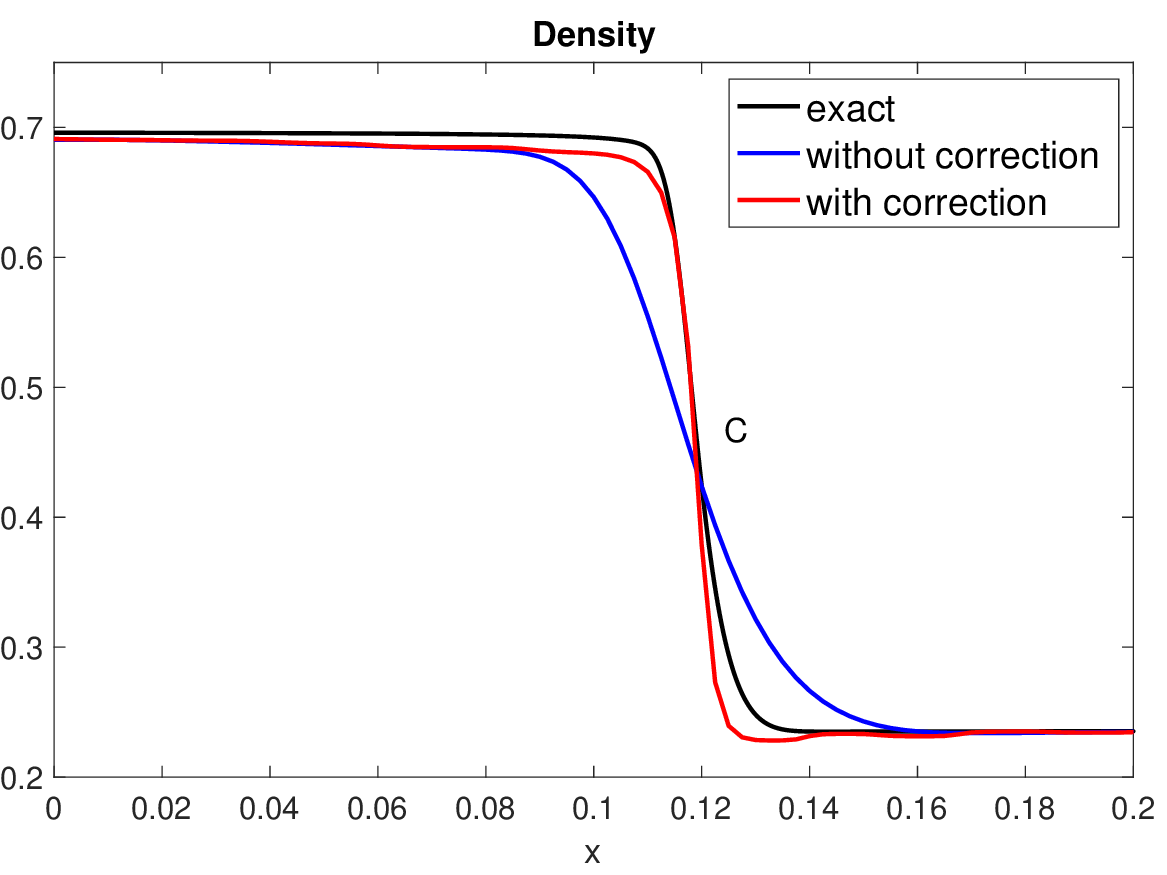}
    }
    \caption{Example~\ref{sec.1DBW} Brio-Wu shock tube problem : Density with 800 grid points, in comparison with the scheme without LDCU projection step.}
    \label{fig.1DBW_com}
\end{figure}

This is a Riemann problem which is taken from \cite{BT} and gives rise to compound waves. The initial data is given by
\begin{equation*}
    (\rho, v_1, v_2, v_3, B_2, B_3, p)^\top = \left\{
    \begin{array}{cc}
       (1.0, 0, 0, 0, 1.0, 0, 1.0)^\top\qquad  &\text{for}\; x<0,\\
       (0.125, 0, 0, 0, -1.0, 0, 0.1)^\top  &\text{for}\; x>0.
    \end{array}
    \right. 
\end{equation*}
with constant $B_1=0.75$ and $\gamma = 2$ and the computational domain is $[-1,1]$. Figure~\ref{fig.1DBW_com} shows the density on 800 grid points at the final time $T=0.2$. The result includes a left-moving fast rarefaction wave (FR), a slow compound wave (SM), a contact discontinuity (C), a right-moving slow shock (SS), and a right-moving fast rarefaction wave (FR). All the waves are well resolved and the second figure shows  that our result resolves the contact discontinuity better than the simulation that does not use the LDCU projection step.

\subsubsection{Dai \& Woodward shock tube problem}
\label{sec.1DDW}
\begin{figure}[htbp]
    \centering   
    \subfigure[]{\includegraphics[height=0.36\linewidth]{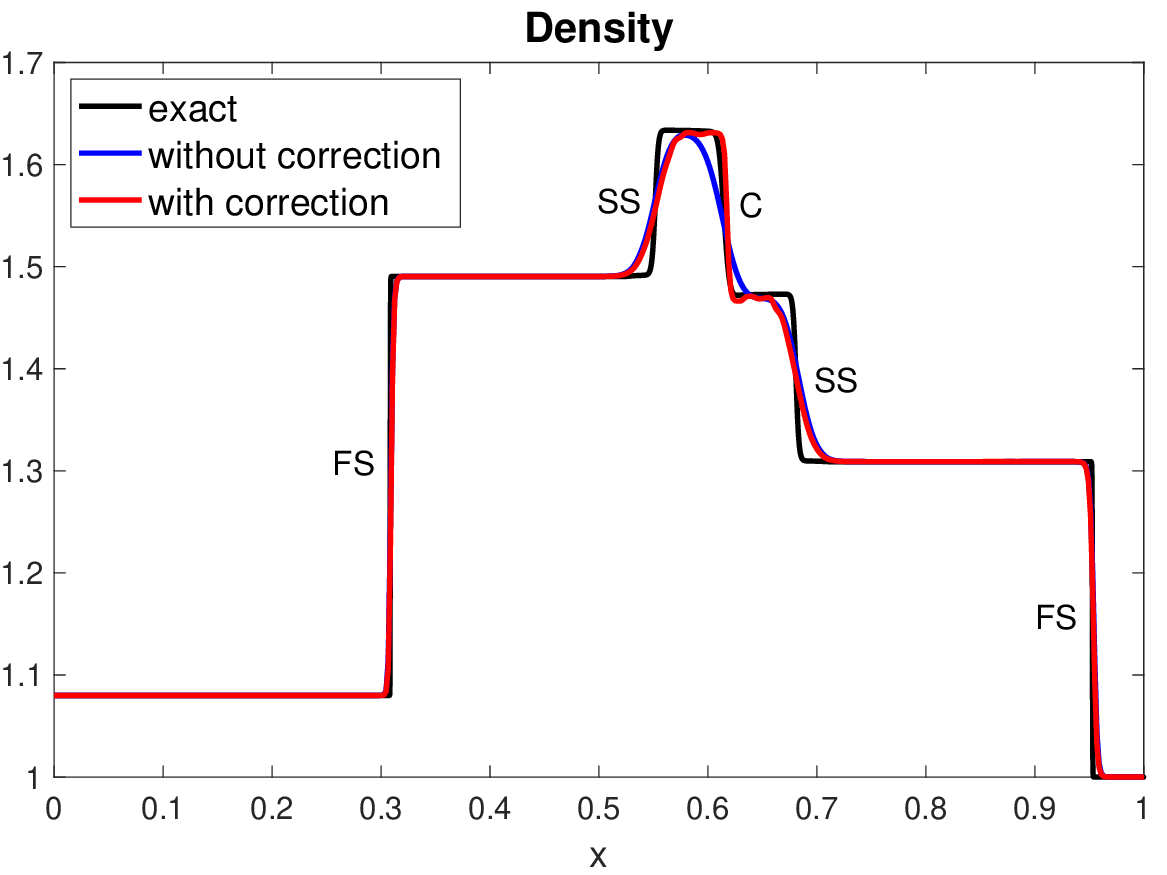} 
    } 
    \subfigure[Zoomed view around contact wave]{\includegraphics[height=0.36\linewidth]{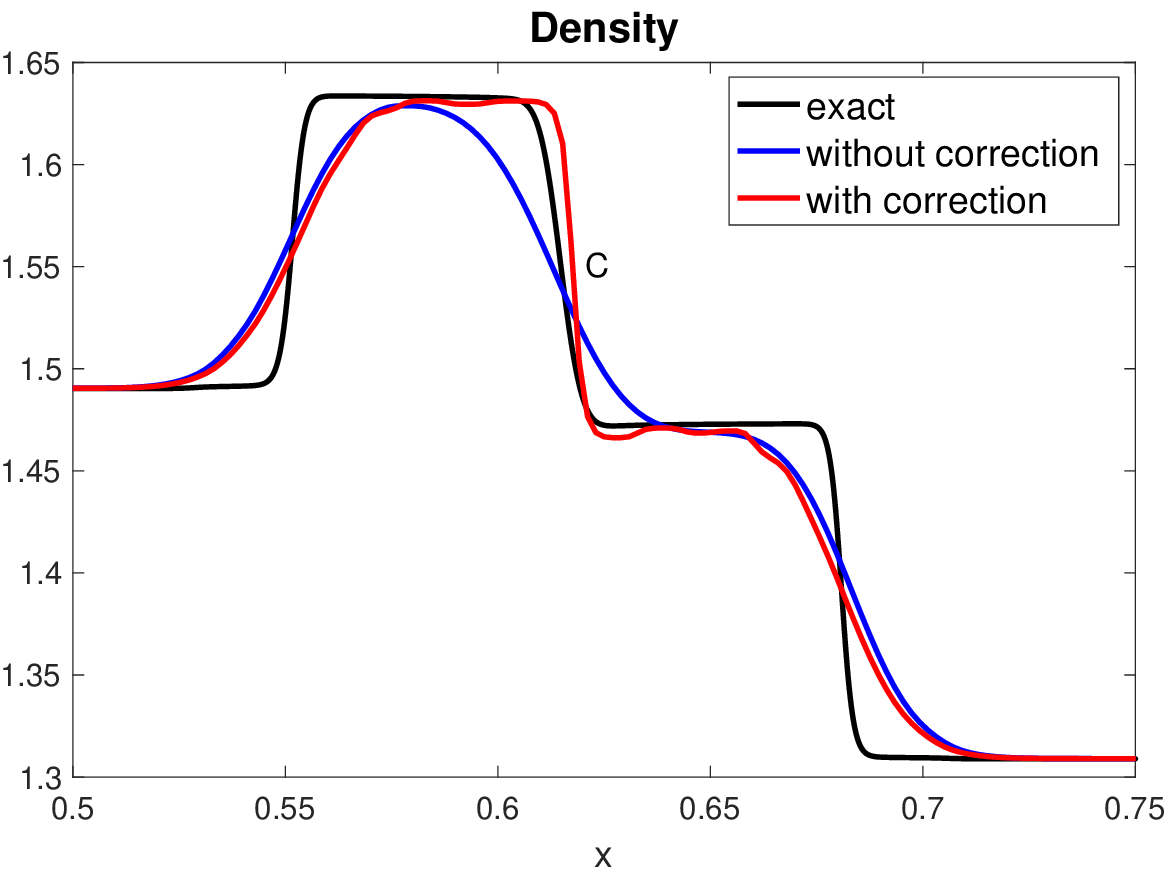}
    }
    \caption{Example~\ref{sec.1DDW} Dai \& Woodward shock tube problem: Density with 512 grid points. Compared with the scheme without LDCU correction.}
    \label{fig.1DDW_com}
\end{figure}

This is a Riemann problem taken from~\cite{1DDW} with initial condition,
\begin{equation*}
    (\rho, v_1, v_2, v_3, B_2, B_3, p)^\top = \left\{
    \begin{array}{cc}
       (1.08, 1.2, 0.01, 0.5, \dfrac{3.6}{\sqrt{4\pi}}, \dfrac{2}{\sqrt{4\pi}}, 0.95)^\top\qquad  &\text{for}\; x<0.5,\\
       (1.0, 0, 0, 0, \dfrac{4}{\sqrt{4\pi}}, \dfrac{2}{\sqrt{4\pi}},1.0)^\top  &\text{for}\; x>0.5.
    \end{array}
    \right. 
\end{equation*}
with $B_1=\frac{2}{\sqrt{4\pi}}$ and $\gamma = \frac{5}{3}$. The computational domain is $[0,1]$. The solution of this problem consists of a pair of fast and slow
left-going magnetosonic waves, (FS) and (SS), a contact discontinuity (C), a pair of fast and slow right-going
magnetosonic waves, and the rotational discontinuity between each pair of
fast and slow magnetosonic waves. The obtained density at the final time $T=0.2$ using 512 grid points is shown in Figure~\ref{fig.1DDW_com}. Again the resolution of the contact discontinuity is much improved.  

\subsubsection{Ryu-Jones problem}
\label{sec.1DRJ}
\begin{figure}[htbp]
    \centering
    \subfigure[]{\includegraphics[height=0.36\linewidth]{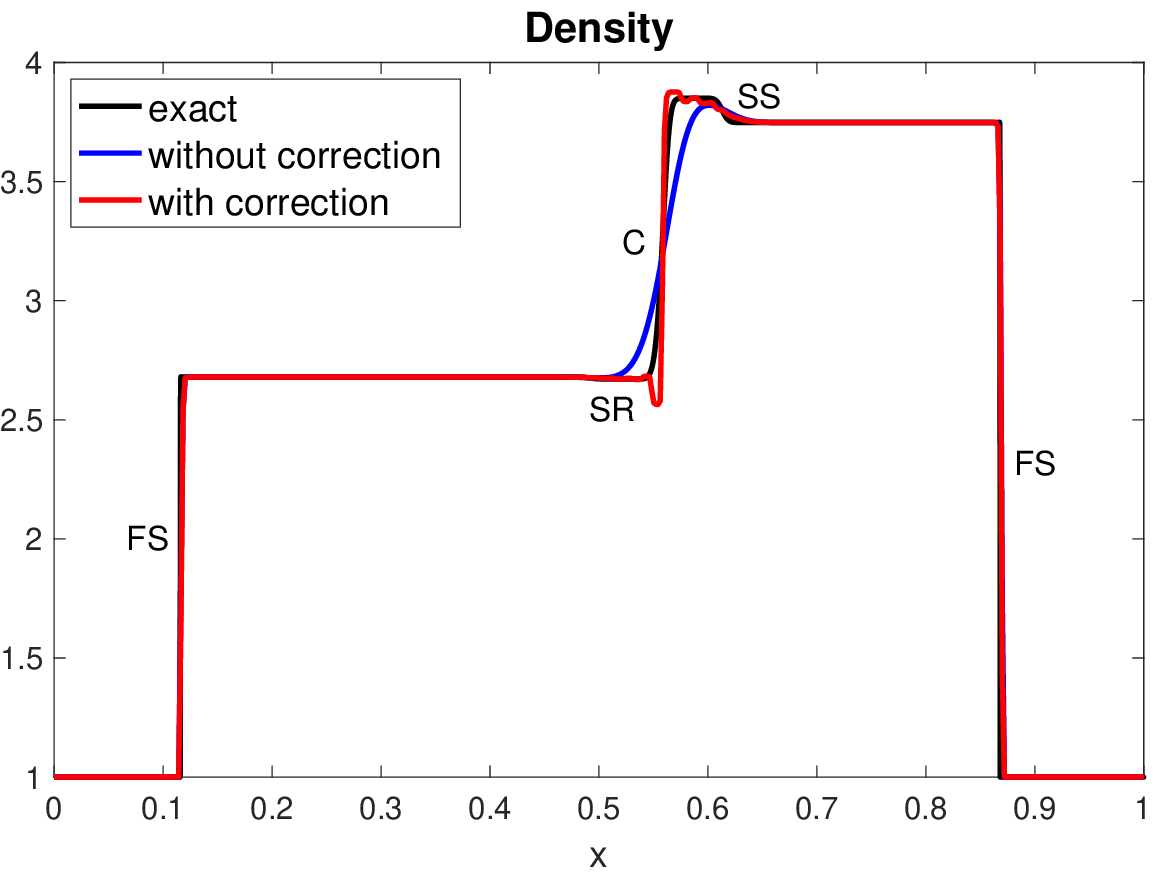}
    }
    \subfigure[Zoomed view around contact wave]{
    \includegraphics[height=0.36\linewidth]{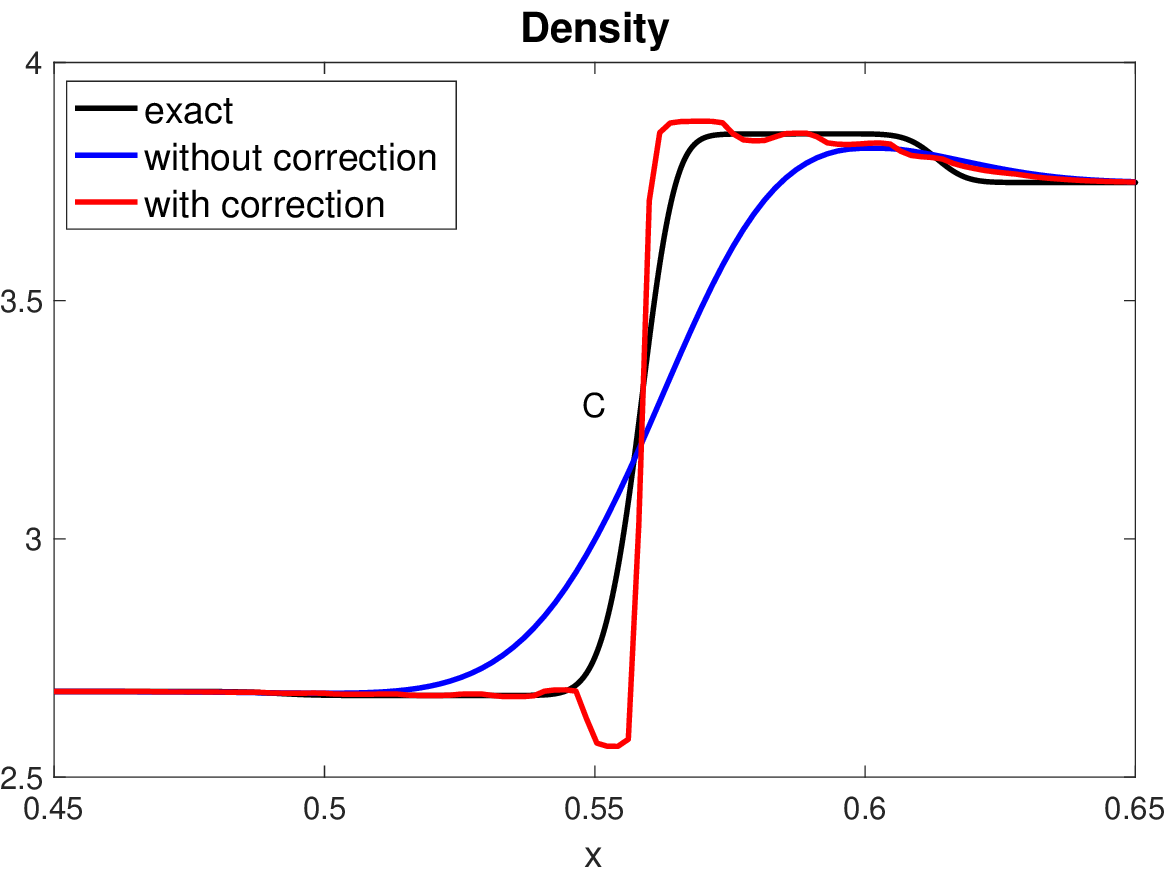}
    }
    \caption{Example~\ref{sec.1DRJ} Ryu-Jones problem: Density with 516 grid points. Compared with the scheme without LDCU correction.}
    \label{fig.1DRJ_com}
\end{figure}

This is a Riemann problem taken from~\cite{1DRJ}; the initial condition is given by
\begin{equation*}
    (\rho, v_1, v_2, v_3, B_2, B_3, p)^\top = \left\{
    \begin{array}{cc}
       (1.0, 10, 0, 0, \dfrac{5}{\sqrt{4\pi}}, 0, 20)^\top\quad  &\text{for}\;\; x\leq0.5,\\
       (1.0, -10, 0, 0,  \dfrac{5}{\sqrt{4\pi}}, 0, 1.0)^\top  &\text{for}\;\; x>0.5.
    \end{array}
    \right. 
\end{equation*}
with $B_1 = \frac{5}{\sqrt{4\pi}}$ and $\gamma = \frac{5}{3}$, and we solve the problem on the domain $[0, 1]$ until a final time $T=0.08$.  Figure~\ref{fig.1DRJ_com} shows the result of density with 516 grid points. The results present that the wave structure of this
problem produce five waves from left to right including a fast shock (FS), a slow rarefaction wave (SR), a contact discontinuity (C), a slow shock (SS), and a fast shock (FS). The result with the correction is close to the reference and captures the discontinuity very well, although there are some undershoots.  

\subsection{2-D experiments}
The 2-D results are obtained using the second order reconstruction in space and the three-stage Runge-Kutta (RK3) method for time evolution. For the ODE, $dQ/dt = C[Q]$, the RK3 method is given by
\begin{align*}
    &Q^{(1)} = Q^n+\Delta tC[Q^n],\\
    &Q^{(2)} = Q^{(1)}+\frac{\Delta t}{4}\left(-3C[Q^n]+C[Q^{(1)}]\right),\\
    &Q^{n+1} = Q^{(2)}+\frac{\Delta t}{12}\left(-C[Q^n]-C[Q^{(1)}]+8C[Q^{(2)}]\right).
\end{align*}
and the same method is applied to Eq.~\eqref{semiQ}, \eqref{semiB1} and \eqref{semiB2}.  We set the CFL number as 0.45, and use the following formula to compute $\Delta t$ in each time step.
\begin{align*}
    \Delta t = \text{CFL} \cdot \min_{j,k}\left(\frac{\Delta x}{\max(a^+_{j+\frac{1}{2},k}, a^-_{j+\frac{1}{2},k})}, \frac{\Delta y}{\max(b^+_{j,k+\frac{1}{2}}, b^-_{j,\frac{1}{2}})}\right).
\end{align*}
A high level view of the 2-D method is described as an algorithm in Section~\ref{sec:algo}.

\subsubsection{Balsara vortex test} 
\label{sec.2Dvortex}
\begin{figure}[htbp]
    \centering
    \subfigure[Exact solution ]{
        \includegraphics[width=0.48\linewidth]{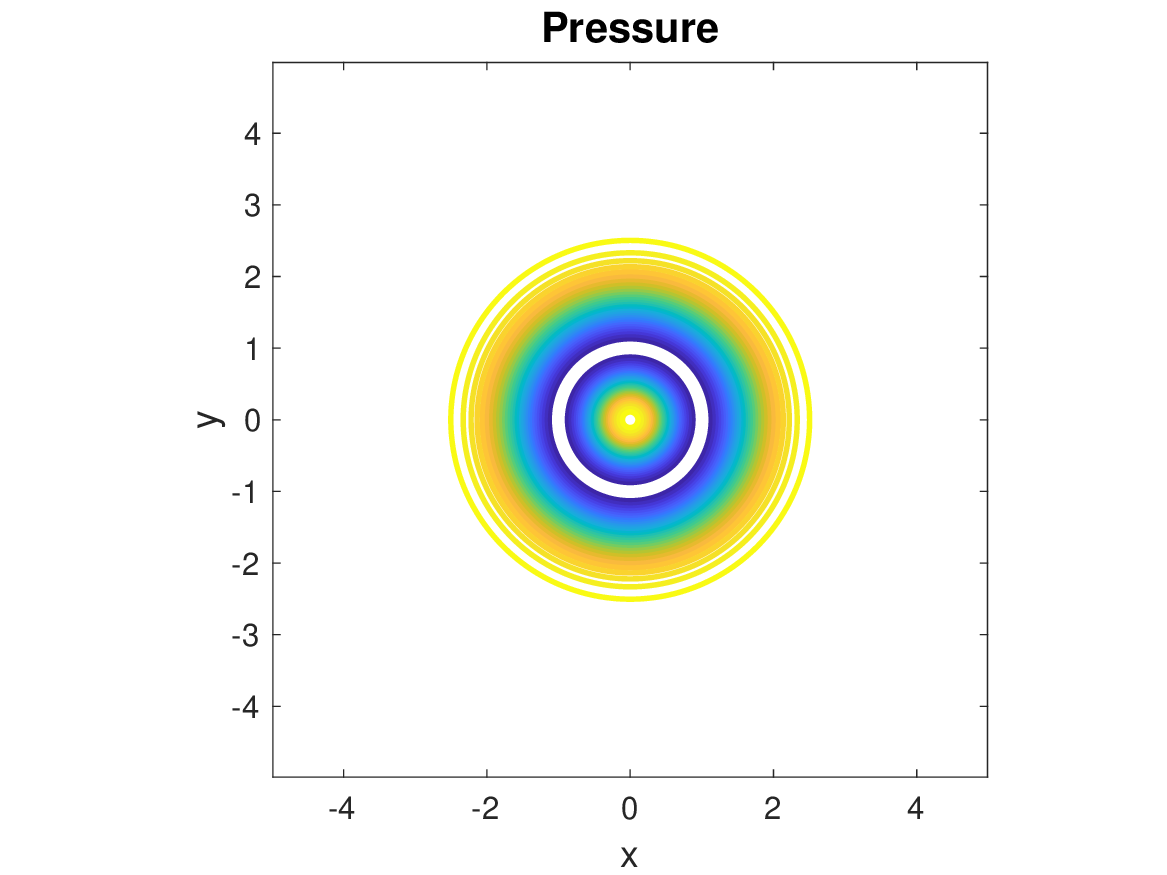}
        } \\  
    \subfigure[$200\times200$]{
        \includegraphics[width=0.48\linewidth]{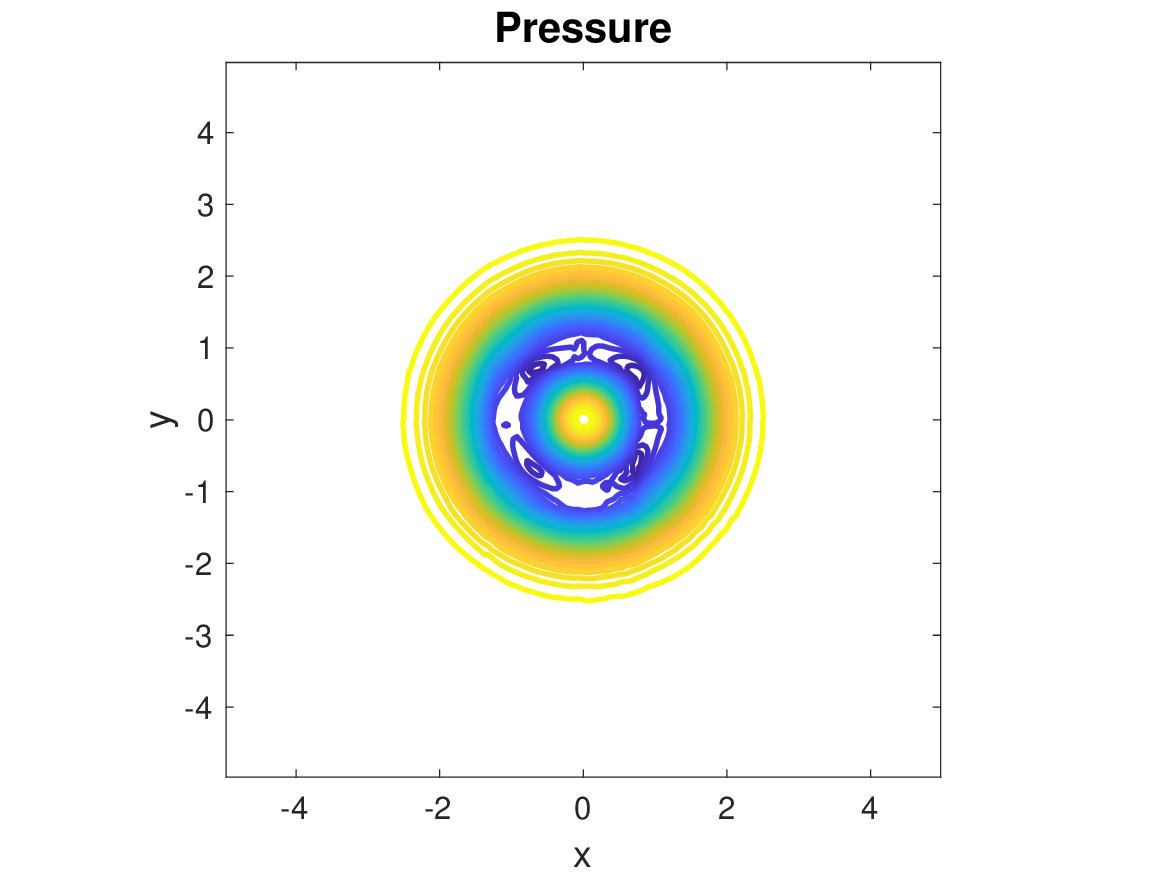}
    }    
    \subfigure[$400\times400$]{
        \includegraphics[width=0.48\linewidth]{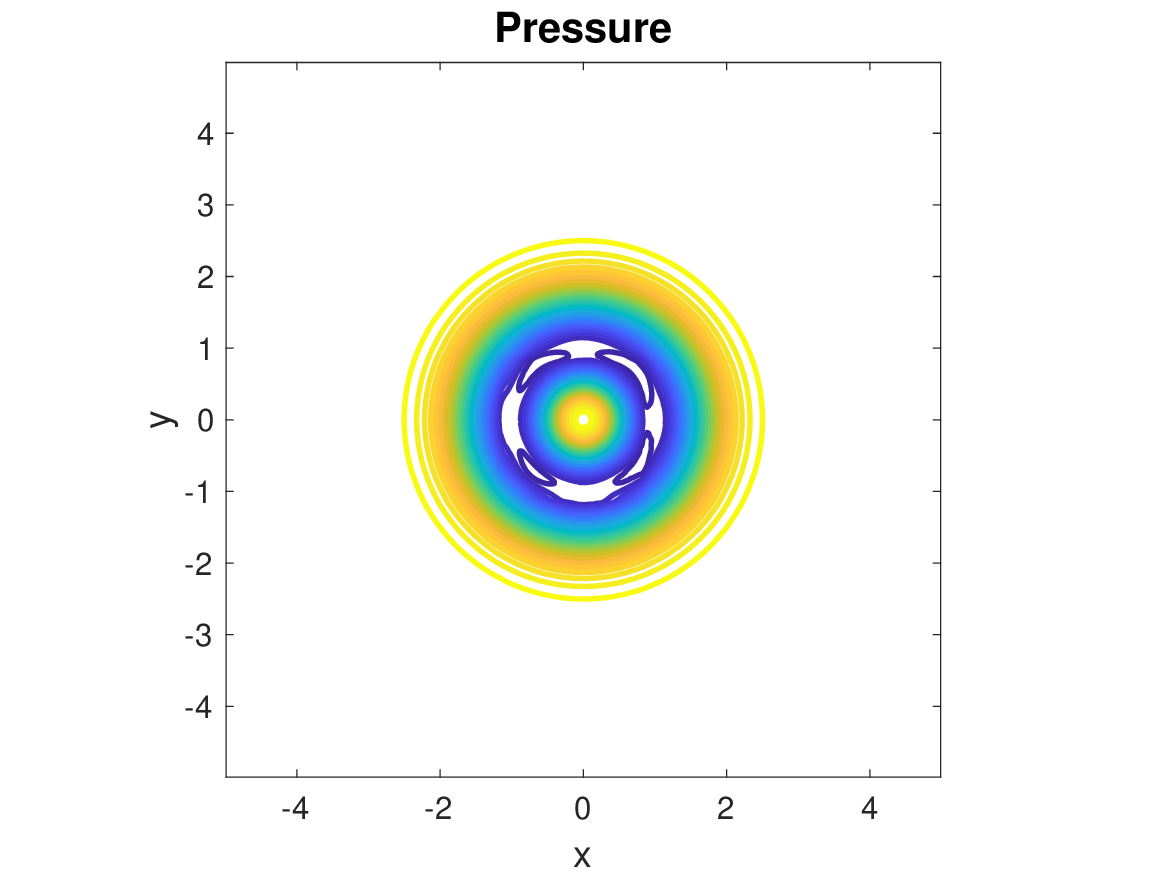}   
        }
    \caption{Example~\ref{sec.2Dvortex} Balsara vortex test: pressure at T=10 with 30 contour lines. The  MC-$\theta$ limiter is used.}
    \label{2Dvotex_WithM}
\end{figure}

\begin{figure}[htbp]
    \centering
    \subfigure[Exact solution]{
        \includegraphics[width=0.48\linewidth]{figs/Result2D/Vortex/p_exact.eps}
        } \\  
    \subfigure[$200\times200$ grid]{
        \includegraphics[width=0.48\linewidth]{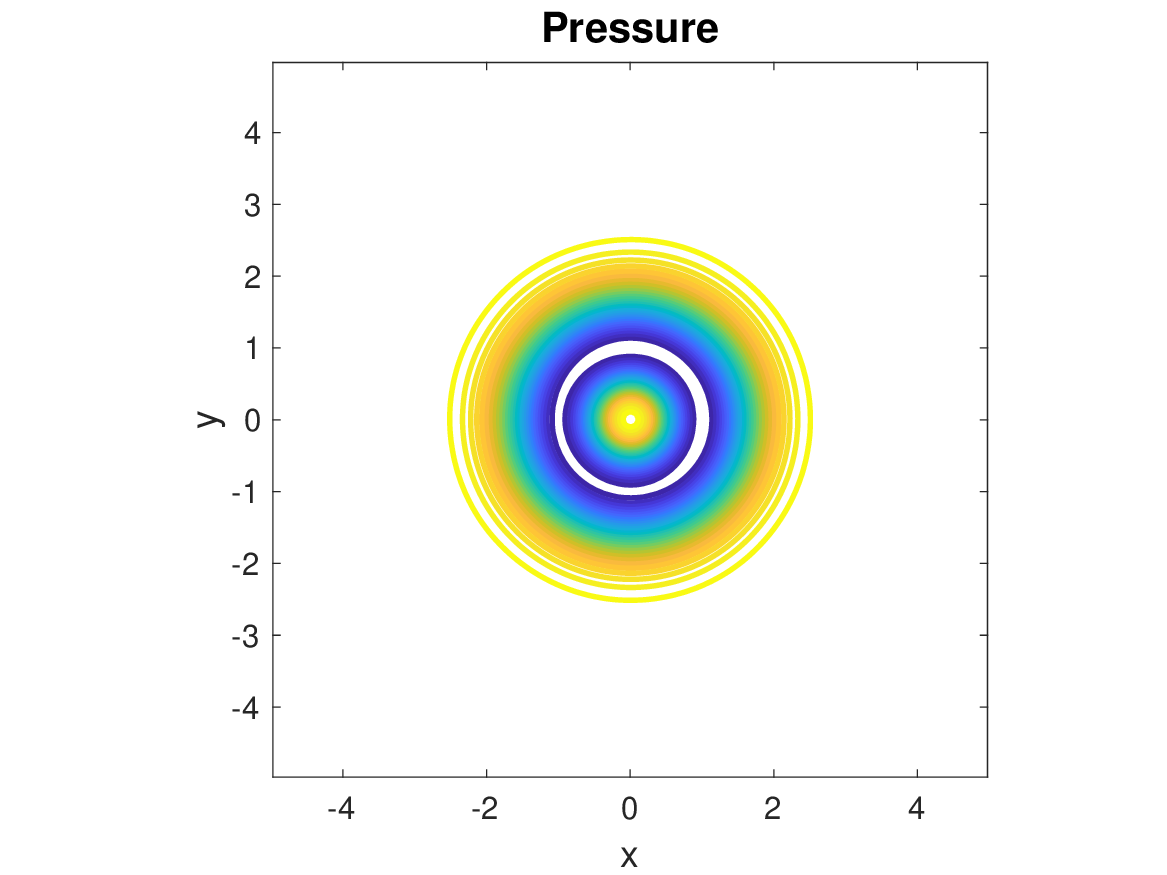}
    }    
    \subfigure[$400\times400$ grid]{
        \includegraphics[width=0.48\linewidth]{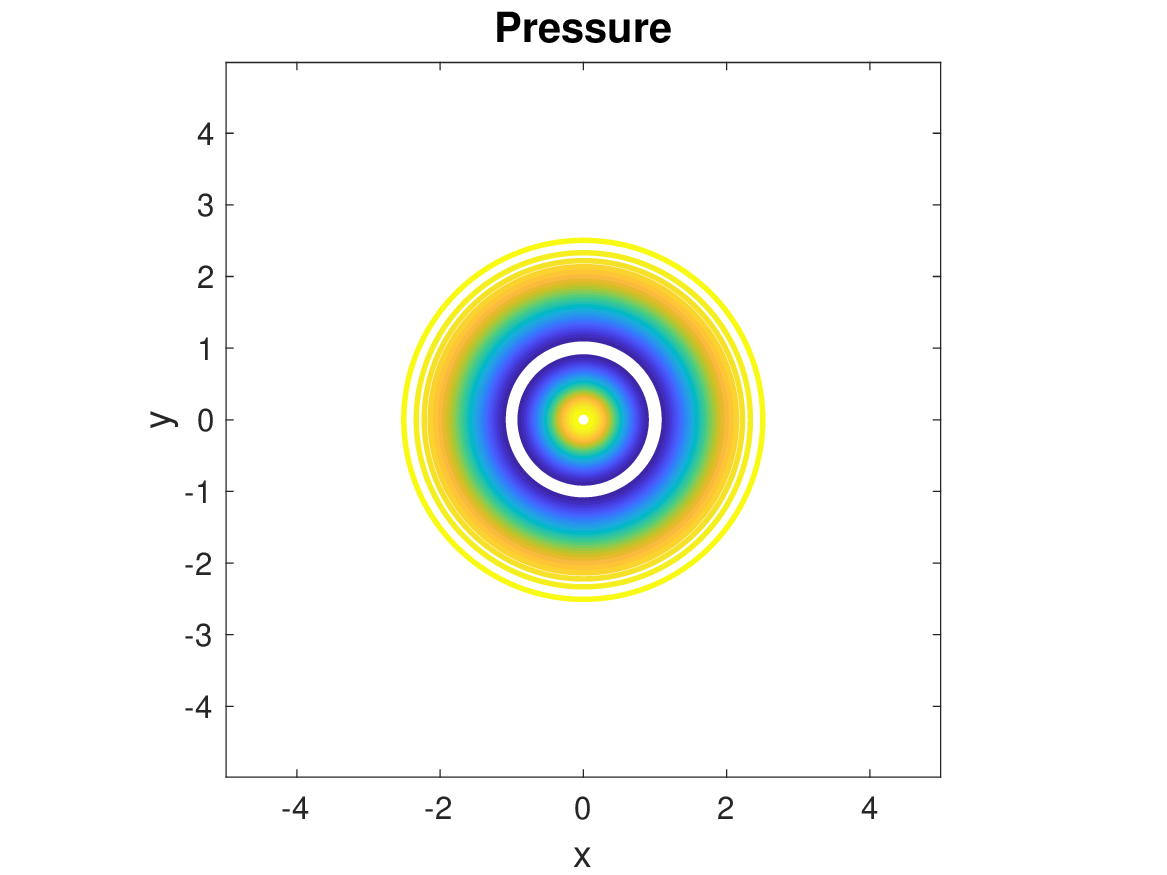}   
        }
    \caption{Example~\ref{sec.2Dvortex} Balsara vortex test: pressure at T=10 with 30 contour lines. No limiter is used.}
    \label{2Dvotex_WithCD}
\end{figure}

In this test, we apply our scheme to the Balsara vortex test from~\cite{ref_Blast} to examine the order of accuracy. We compute in the domain $[-5,5] \times [-5,5]$ with periodic boundary condition. Following the setup in \cite{ref_Balsara_setup}, the initial data is the background flow $(\rho,v,B_1,B_2,B_3,p) = (1,1,1,0,0,0,0,1)$ with the perturbation
\begin{align*}
    (\delta v_1,\delta v_2) = \xi \exp(0.5(1-r^2))(-y,x),\\
    (\delta B_1,\delta B_2) = \mu \exp(0.5(1-r^2))(-y,x), \\
    \delta p= 0.5(\mu^2(1-r^2)-\xi^2) \exp(1-r^2),
\end{align*}
where $r=\sqrt{x^2+y^2}$ $\zeta=\frac{1}{2\pi}$, and $\mu=\frac{1}{2\pi}$. Figure \ref{2Dvotex_WithM} and Table \ref{Votex_err_WithM} show the $l^1$-error and the  convergence rate at the  final time $T=10$ by the proposed scheme with the reconstruction using the MC-$\theta$ limiter, while   the result in Figure \ref{2Dvotex_WithCD} and Table \ref{Vortex_err_withCD} using the reconstruction with central difference, i.e., without limiter.  Both of the results show that our scheme can achieve close to second-order accuracy.

\begin{table}
\caption{Example~\ref{sec.2Dvortex} Balsara vortex test: $l^1$-error and convergence rate with MC-$\theta$ limiter. 
}\label{Votex_err_WithM}
\centering
\begin{tabular}{lllll}
\toprule
& $\rho$-Error & $\rho$-Order & P-Error & P-Order\\  
\midrule
 50x50 & 5.1x$10^{-2}$ & - & 4.0x$10^{-2}$ & - \\
 100x100 & 1.7x$10^{-2}$ & 1.58 & 1.0x$10^{-2}$ & 2.00 \\
 200x200 & 4.0x$10^{-3}$ & 2.08 & 2.6x$10^{-3}$ & 1.97 \\
 400x400 & 8.8x$10^{-4}$ & 2.21 & 6.6x$10^{-4}$ & 1.95 \\
\bottomrule
\end{tabular}
\end{table}

\begin{table}
\caption{Example~\ref{sec.2Dvortex} Balsara vortex test: $l^1$-error and convergence rate without limiter.}
\label{Vortex_err_withCD}
\centering
\begin{tabular} {lllll}
\toprule
 & $\rho$-Error & $\rho$-Order & P-Error & P-Order\\ 
 \midrule
 50x50 & 7.8x$10^{-3}$ & - & 2.5x$10^{-2}$ & - \\
100x100 & 2.2x$10^{-3}$ & 1.83 & 5.8x$10^{-3}$ & 2.13 \\ 
200x200 & 5.4x$10^{-4}$ & 2.03 & 1.3x$10^{-3}$ & 2.12 \\ 
400x400 & 1.3x$10^{-4}$ & 2.08 & 3.2x$10^{-4}$ & 2.06 \\ 
\bottomrule
\end{tabular}
\end{table}

\subsubsection{Smooth sine wave}
\label{sec.2Dsine}
We adopt the smooth sine wave test from \cite{SmoothP} for the second order accuracy test.
Consider the solution whose exact solution is given by
\begin{equation*}
    (\rho, v_1, v_2,v_3,B_1,B_2,B_3,p)^\top = (1+0.99\sin(2\pi(x+y-2t)), 1, 1, 0, 0.1, 0.1, 0, 1)^\top.
\end{equation*}
and we compute up to the final time $T=0.1$. The computational domain is $[0,1]\times[0,1]$ with the periodic boundary condition. Table \ref{table. SineWave} shows the $l^1$-error and the convergence rate of the density, and demonstrate that our scheme is second-order accuracy.

\begin{table}
\caption{Example~\ref{sec.2Dsine} Smooth sine wave: Density $l^1$-error and order at $T=0.1$ with minmod limiter.} \label{table. SineWave}
\centering
\begin{tabular} {lll}
\toprule
 & Error & Order\\ 
 \midrule
50x50 & 3.4x$10^{-1}$ & - \\ 
100x100 & 8.3x$10^{-2}$ & 2.05  \\ 
200x200 & 1.9x$10^{-2}$ & 2.13  \\ 
400x400 & 4.5x$10^{-3}$ & 2.07  \\ 
\bottomrule
\end{tabular}
\end{table}

\subsubsection{Brio-Wu shock tube problem}
\label{sec.2DBW}
\begin{figure}[htbp]
    \centering        \includegraphics[width=0.5\linewidth]{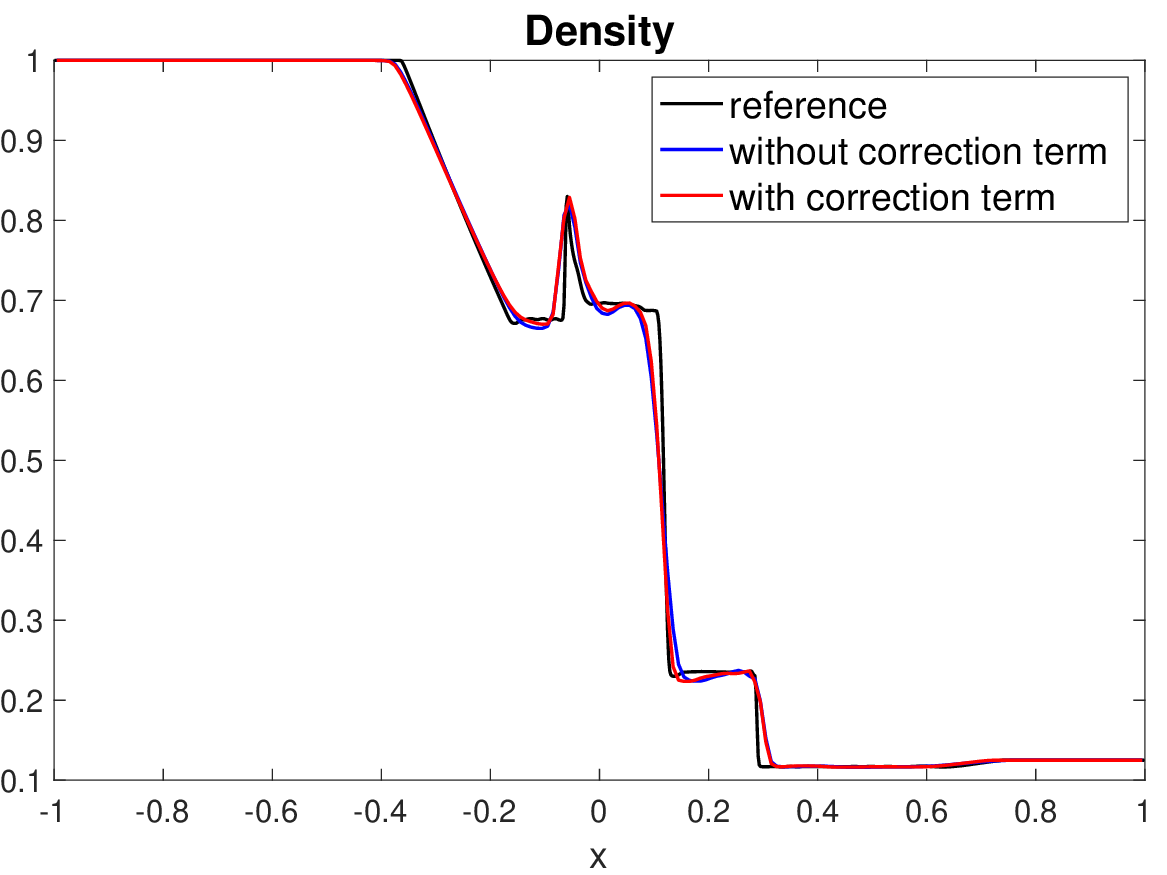}
    \caption{Example~\ref{sec.2DBW} Brio-Wu shock tube problem: The cross-section of the density on $200\times200$ grid points.}
    \label{fig.2DBW}
\end{figure}

In the first test, we use the same initial data as in 1-D case in Section \ref{sec.1DBW},
\begin{equation*}
    (\rho, v_1, v_2, v_3, B_1, B_2, B_3, p)^\top = \left\{
    \begin{array}{cc}
       (1.0, 0, 0, 0, 0.75, 1.0, 0, 1.0)^\top\qquad  &\text{for}\; x<0,\\
       (0.125, 0, 0, 0, 0.75, -1.0, 0, 0.1)^\top  &\text{for}\; x>0,
    \end{array}
    \right. 
\end{equation*}
with $\gamma =2$, and compute the solution in the domain $[-1,1]\times[-1,1]$.  Figure \ref{fig.2DBW} shows the result of the density on $200\times200$ grid points at the final time $T=0.2$, in comparison to the scheme without the correction term K on $200\times200$ grid points and the finer grid points $1000\times1000$ as the reference solution.  The result with correction term performs better at the contact than the one without correction term. 

\subsubsection{Orszag-Tang MHD turbulence problem}
\label{sec.2DOST}
\begin{figure}[htbp]
    \centering
    \subfigure[ $200\times200$ grid]{
        \includegraphics[width=0.48\linewidth]{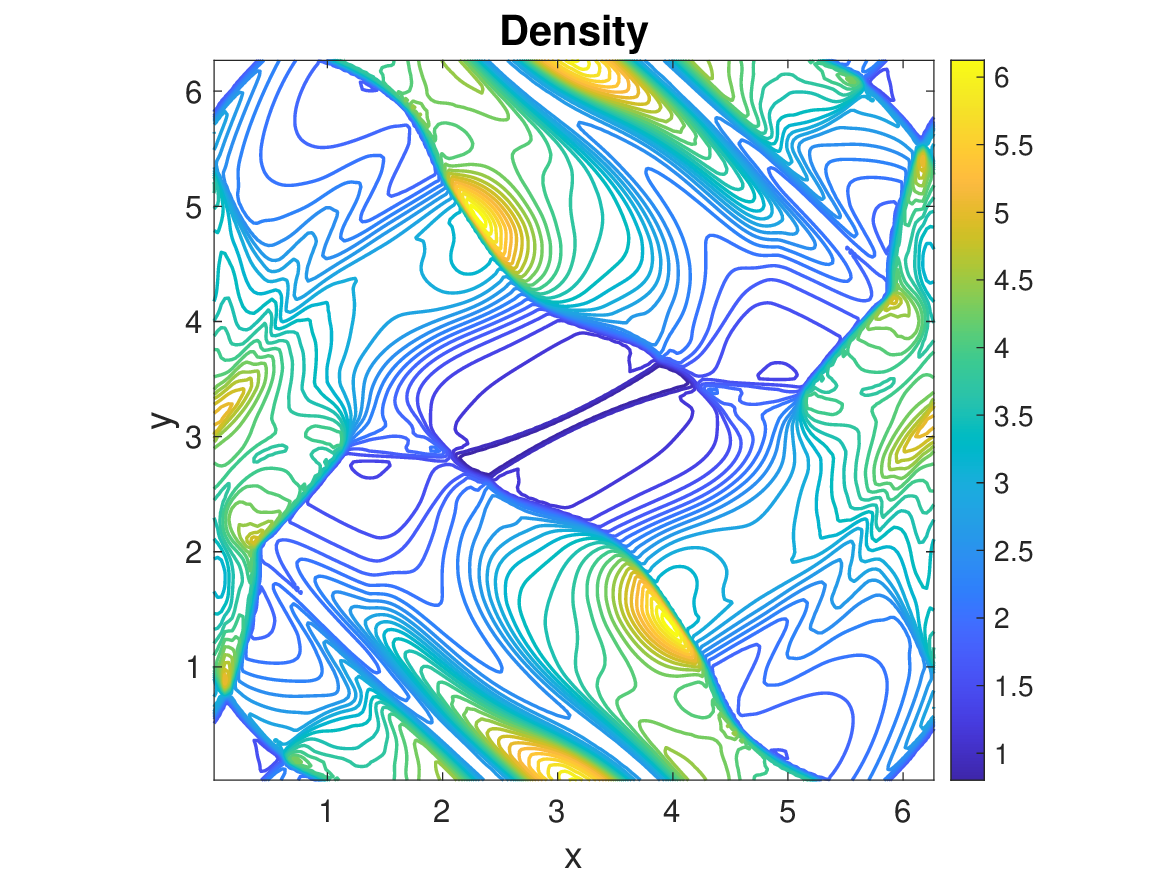}
        } 
    \subfigure[ $400\times400$ grid]{
        \includegraphics[width=0.48\linewidth]{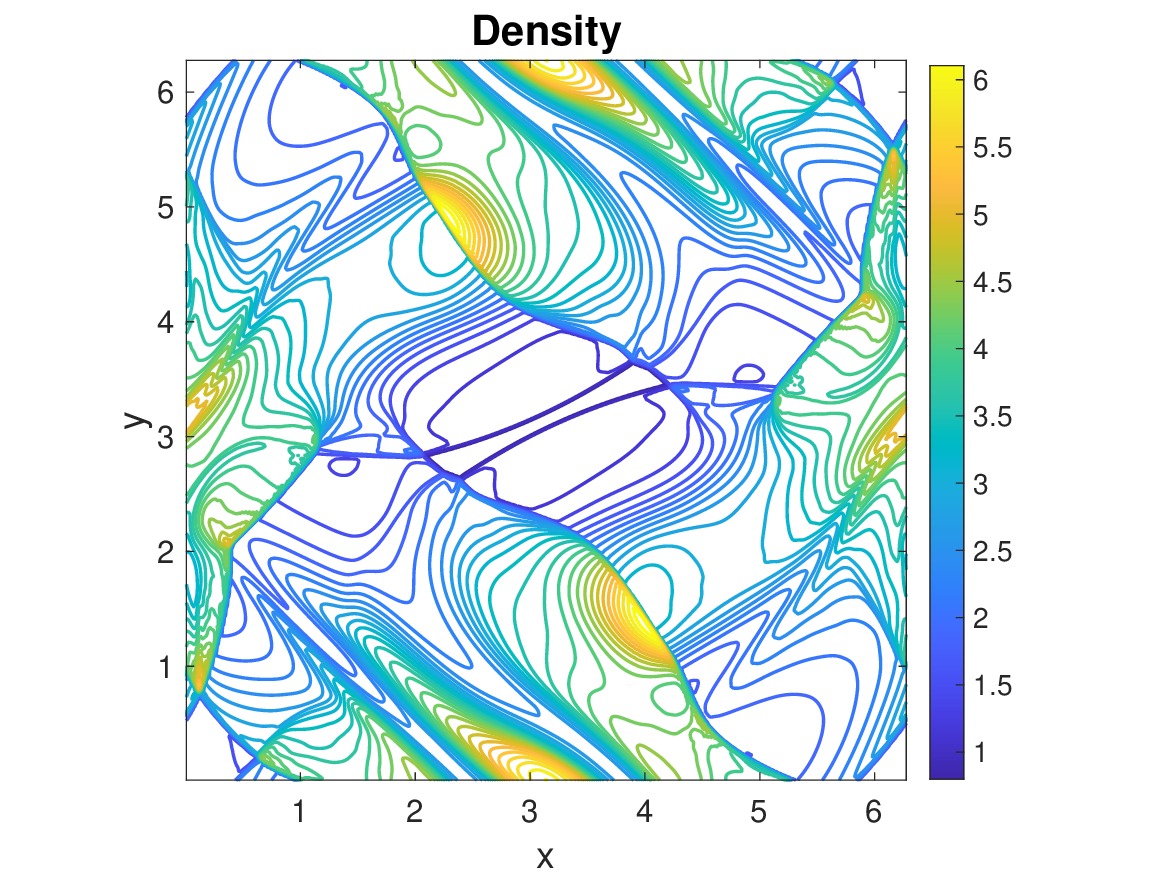}
    }   
    \caption{Example~\ref{sec.2DOST} Orszag-Tang test: Density at $T=2$ with 30 contour lines.}
    \label{fig.OST_T2}
\end{figure}

\begin{figure}[htbp]
    \centering
    \subfigure[Result without correction term on $200\times200$.]{
        \includegraphics[width=0.48\linewidth]{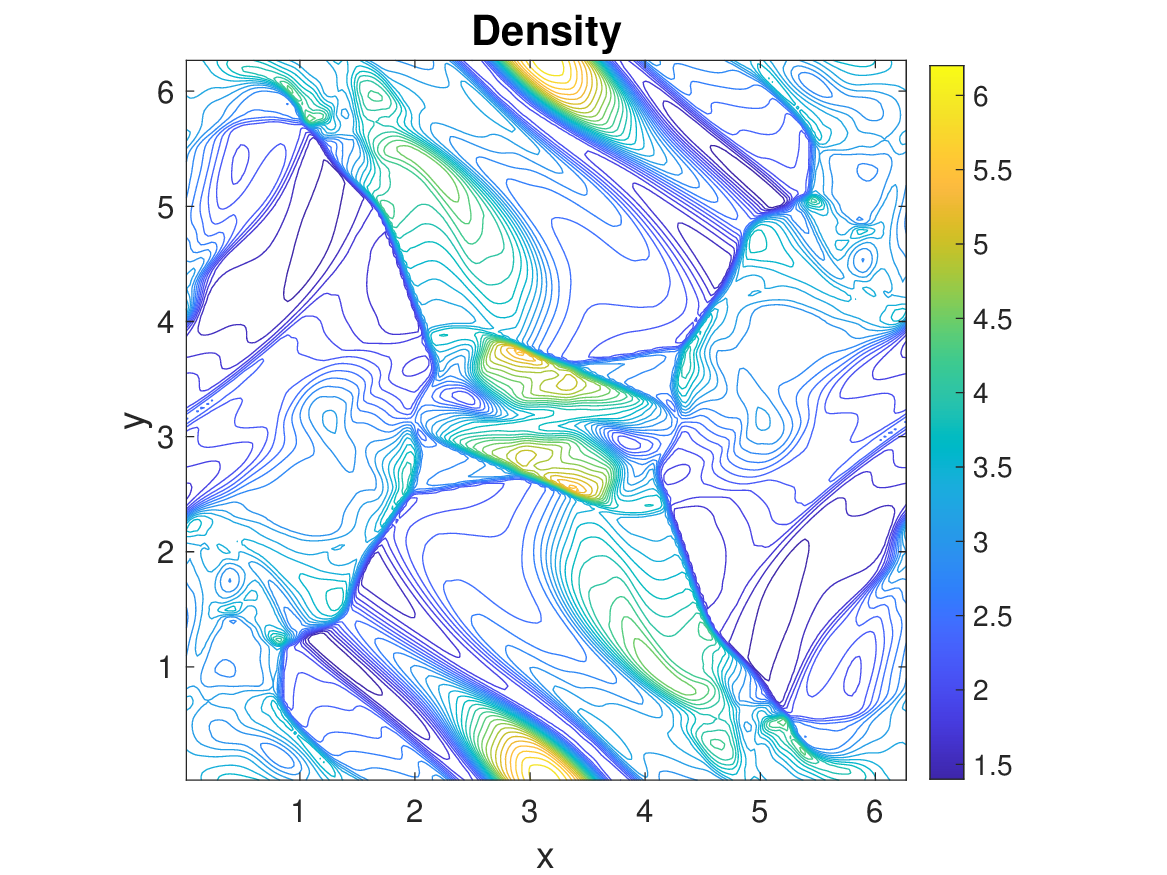}
        } 
    \subfigure[Result with correction term on $200\times200$.]{
        \includegraphics[width=0.48\linewidth]{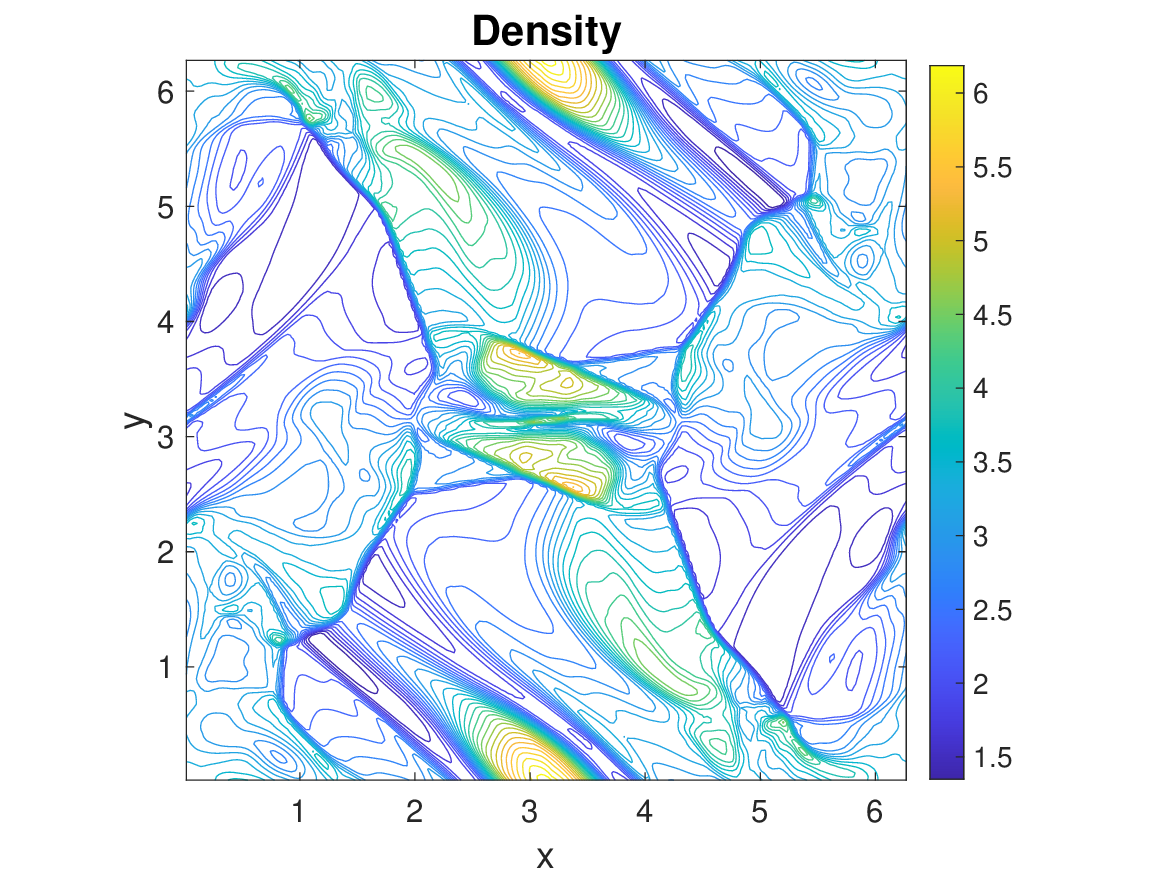}
    }  \\
    \subfigure[Cross-section along $x=\pi$]{
        \includegraphics[width=0.75\linewidth]{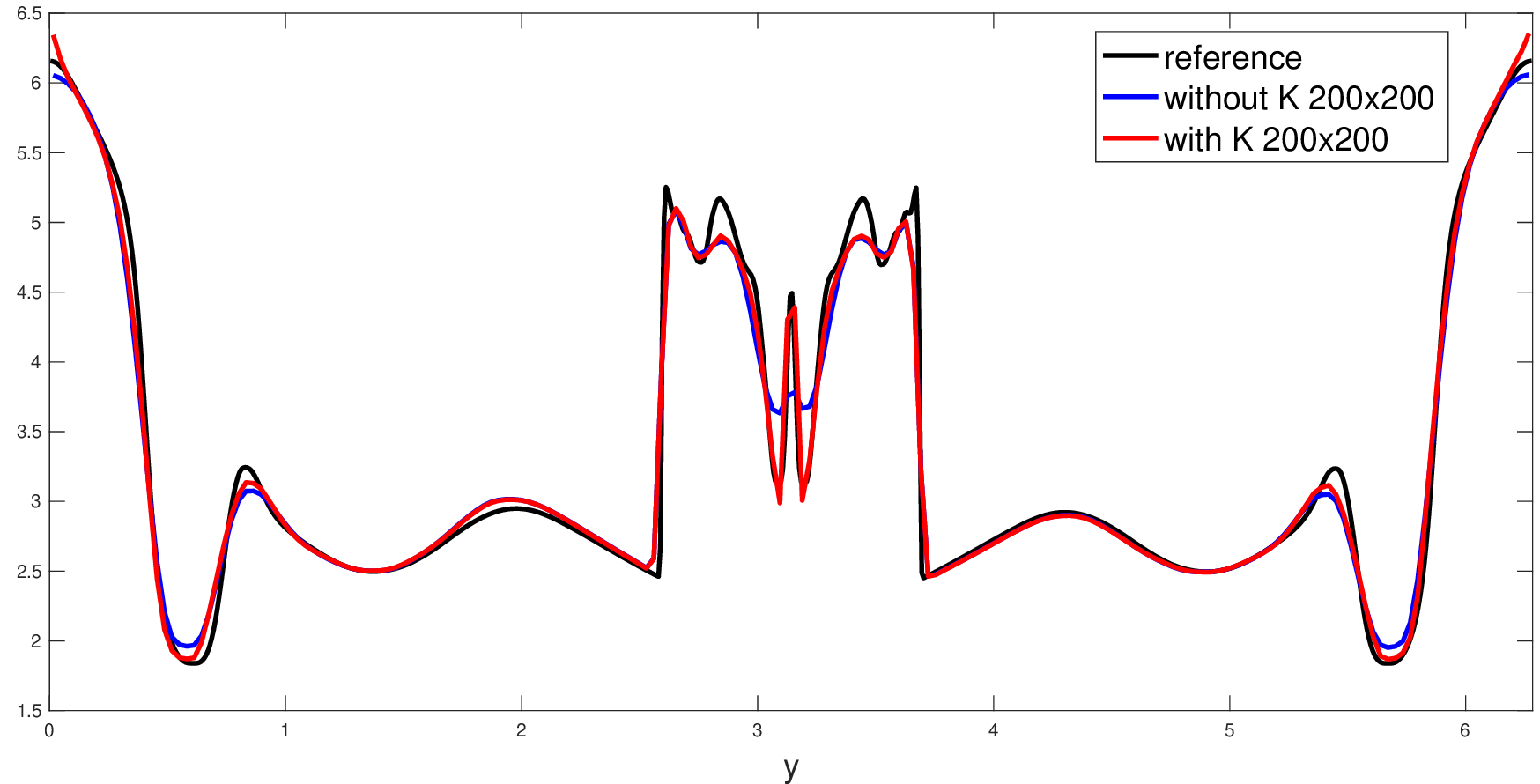}
    }\\
    \subfigure[Zoom in of (c) in $y={[2.4,3.8]}$  ]{
    \includegraphics[width=0.75\linewidth]{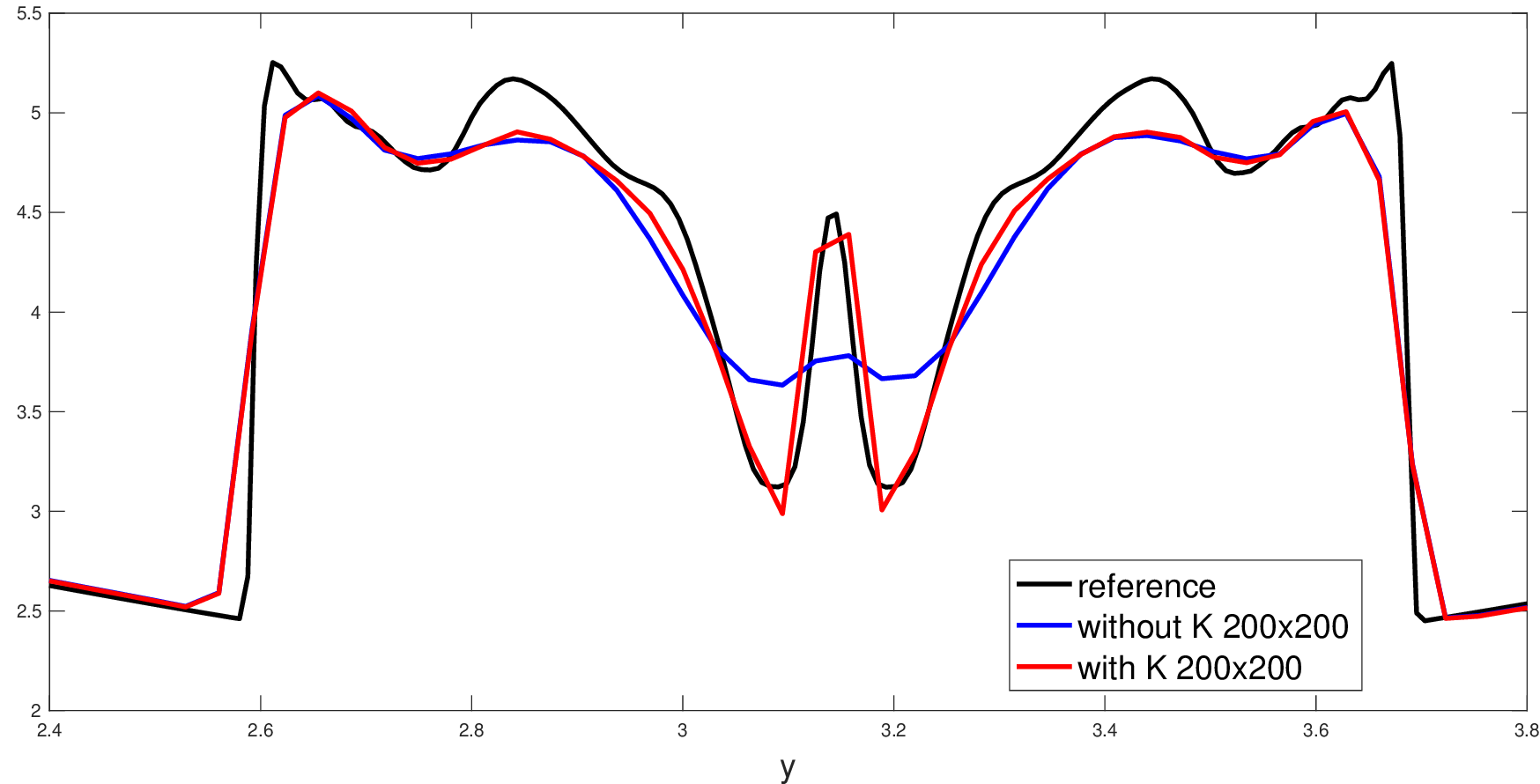}
    }
    \caption{Example \ref{sec.2DOST} Orszag-Tang problem: Comparison of density at $T=3$. Reference solution is obtained by the scheme without K on $800\times800$ grid points. (a) and (b) are drawn with 30 contour lines. }
    \label{fig.OST_T3}
\end{figure}

\begin{figure}[htbp]
    \centering
    \subfigure[$200\times200$ grid]{
        \includegraphics[width=0.48\linewidth]{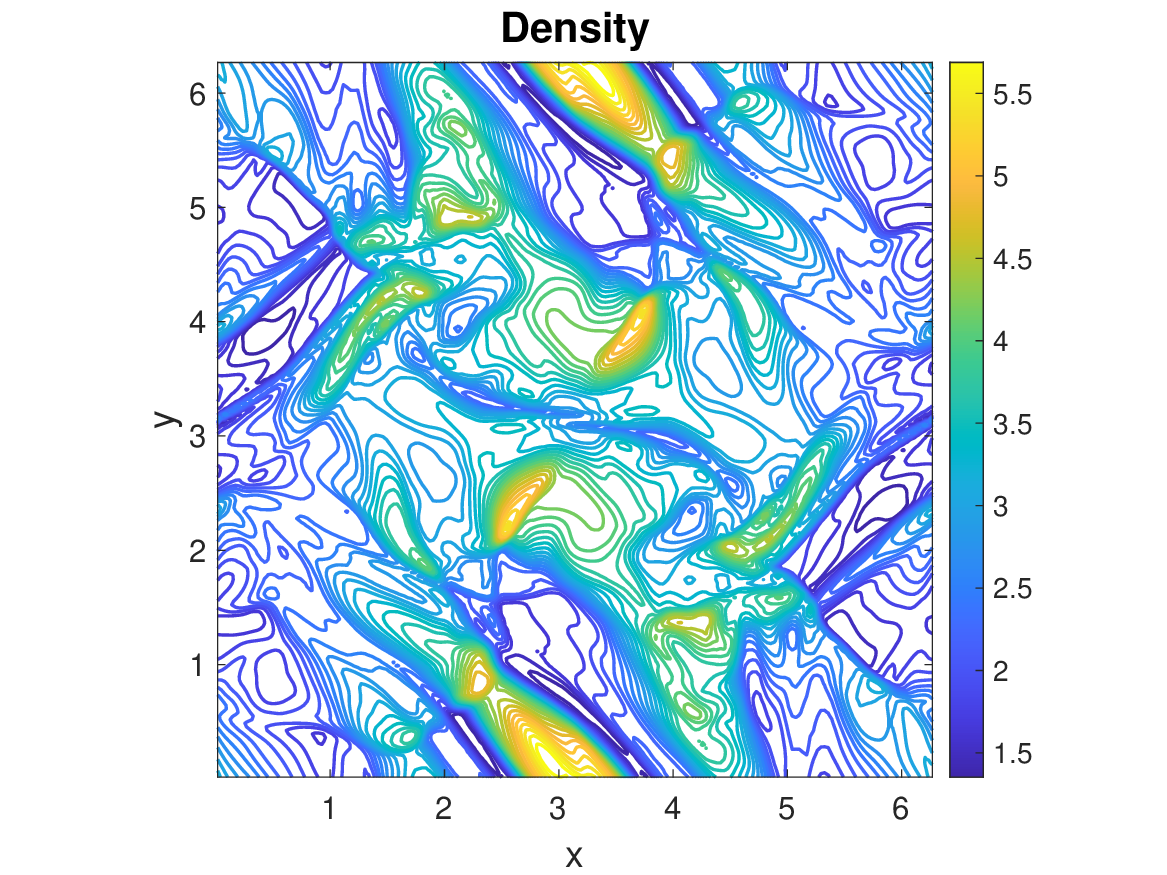}
        } 
    \subfigure[$400\times400$ grid]{
        \includegraphics[width=0.48\linewidth]{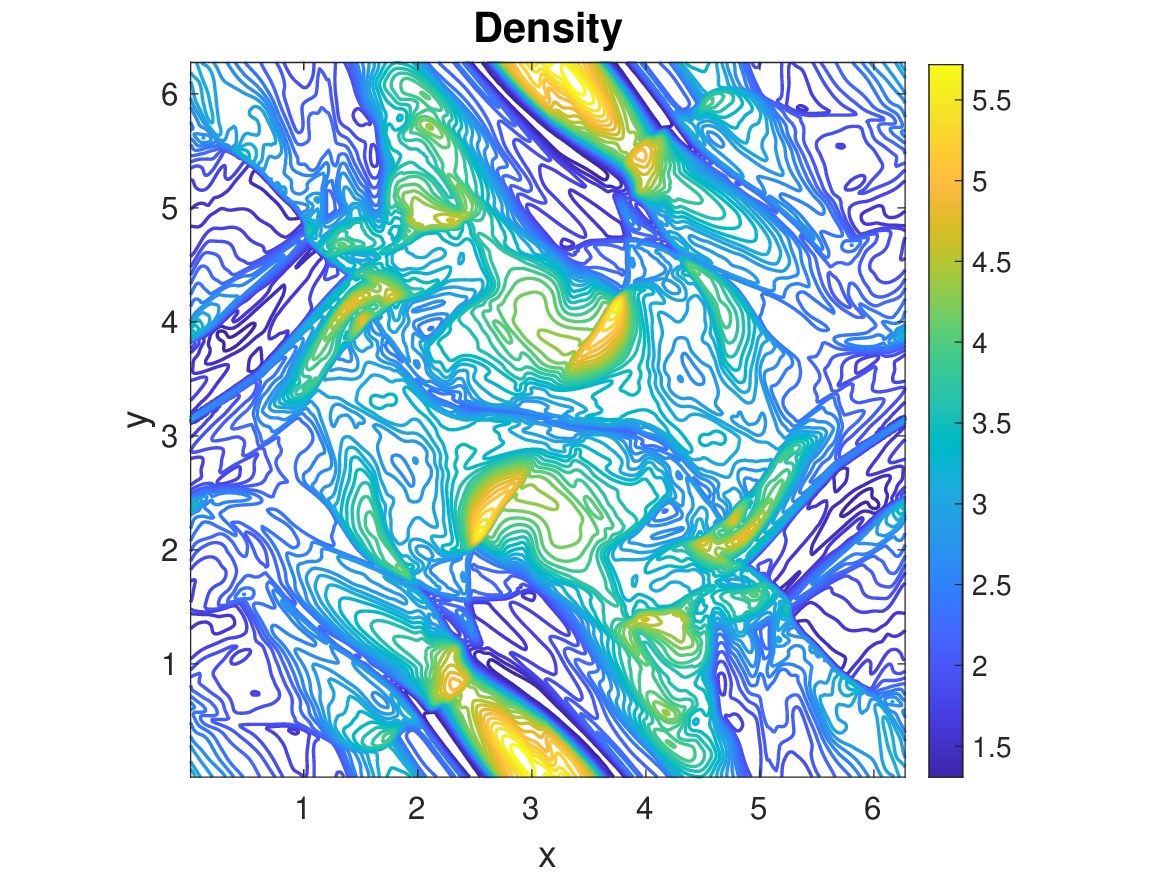}
    }  \\
    \subfigure[Semilogy plot of maximum relative divergence of B versus time]{
        \includegraphics[width=0.8\linewidth]{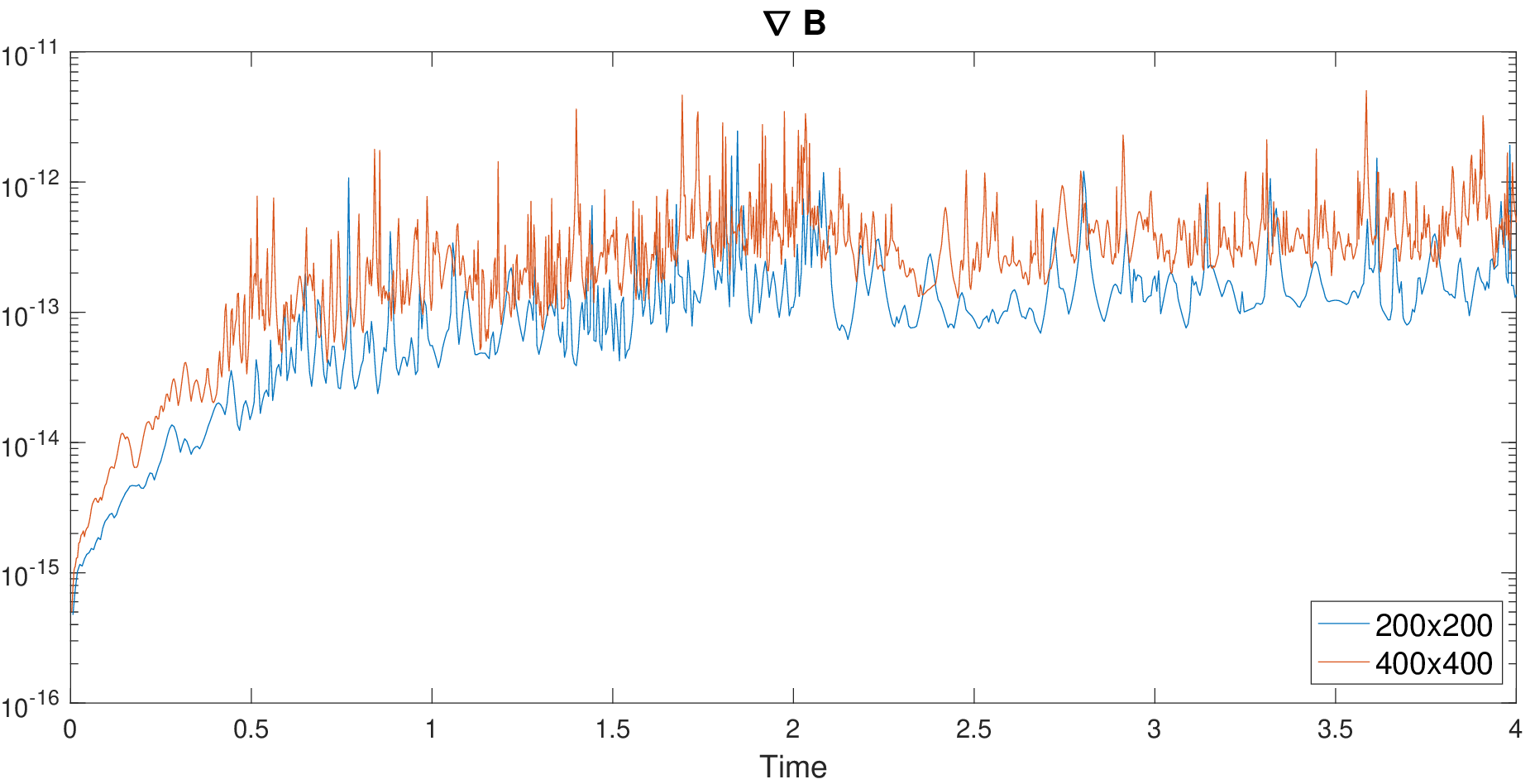}
        }
    \caption{Example \ref{sec.2DOST} Orszag-Tang problem: (a), (b) density with 30 contour lines at $T=4$, (c) evolution of divergence of \textbf{B} with time.}
    \label{fig.OST_T4}
\end{figure}

This test case is commonly used to test the ability of schemes to resolve multiple waves and eventual production of MHD turbulence. We consider the initial data from \cite{BT},
\begin{align*}
    &\rho(x,y,0) = \gamma^2,\; v_1(x,y,0)=-\sin y,\; v_2(x,y,0) = \sin x, \\
    &p(x,y,0) = \gamma,\; B_1(x,y,0)=-\sin y,\; B_2(x,y,0) = \sin 2x,
\end{align*}
where $\gamma = \frac{5}{3}$. The test case is computed in $[0,2\pi]\times[0,2\pi]$ with periodic boundary conditions. Figure~\ref{fig.OST_T2} shows the results of density on $200\times200$ and $400\times400$ grid points at time $T=2$.  We compare the result of density on $200\times200$ grid points at $T=3$ using the proposed scheme with the scheme without correction term K, and plot the cross-section along $x=\pi$ in Figure \ref{fig.OST_T3}. We can observe that the LDCU scheme results in better solutions as they are closer to the reference solution than the uncorrected scheme. In Figure \ref{fig.OST_T4}, we present the results of density on $200\times200$ and $400\times400$ grid points at $T=4$, and plot the time evolution of the 
\[
\textrm{maximum relative divergence} = \max_{j,k}\frac{|\nabla\cdot \mathbf{B}|}{|\mathbf{B}|}\Delta x
\]

where the divergence of \textbf{B} is obtained by
\begin{equation*}
    (\nabla\cdot \textbf{B})_{j,k} = \frac{(B_1)_{j+\frac{1}{2},k}-(B_1)_{j-\frac{1}{2},k}}{\Delta x} + \frac{(B_2)_{j,k+\frac{1}{2}}-(B_2)_{j,k-\frac{1}{2}}}{\Delta y}.
\end{equation*}
Notice that the deviation of $\mathbf B$ from $\nabla \cdot \mathbf{B} = 0$ stays very small over time and does not increase in time.


\subsubsection{Rotor Test}
\label{sec.2Drotor}
We consider the second rotor problem from \cite{Toth}. The initial data of this test is given by
\begin{equation*}
    \rho= \left\{
    \begin{array}{cc}
        10 \qquad\qquad r\leq r_0, \\
        1+9f(r)\quad  r_0< r <r_1, \\
        1\qquad\;\;\qquad  r\geq r_1,
    \end{array}
    \right.
\end{equation*}
\begin{equation*}
    (v_1, v_2, v_3)= \left\{
    \begin{array}{cc}
        \dfrac{v_0}{r_0}(-(y-\dfrac{1}{2}), x-\dfrac{1}{2}, 0)\quad &r\leq r_0, \\
        \dfrac{f(r)v_0}{r_0}(-(y-\dfrac{1}{2}), x-\dfrac{1}{2}, 0)\quad & r_0< r <r_1, \\
        (0,0, 0)\quad & r\geq r_1,
    \end{array}
    \right.
\end{equation*}
\begin{equation*}
    (B_2, B_3) = (0,0),
\end{equation*}
where $r_0=0.1$, $r_1=0.115$, $r=\sqrt{(x-\frac{1}{2})^2+(y-\frac{1}{2})^2}$, and $f(r)=\frac{r_1-r}{r_1-r_0}$, and $B_1=\frac{2.5}{\sqrt{4\pi}}$, $p=0.5$, $v_0=1$.  The final time is $T=0.295$.

\begin{figure}[htbp]
    \centering
    \subfigure[Without K on $200\times200$]{
    \includegraphics[width=0.48\linewidth]{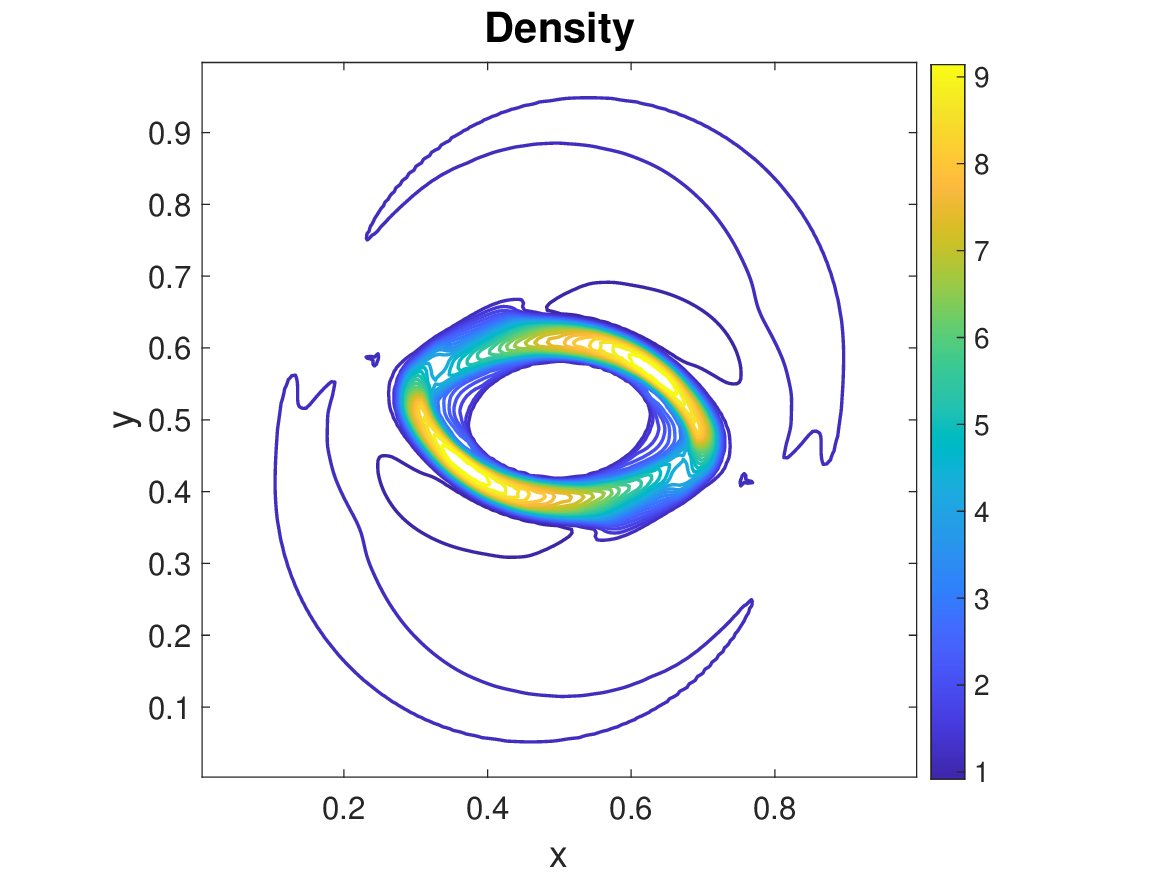}
    }
    \subfigure[With K on $200\times200$]{
    \includegraphics[width=0.48\linewidth]{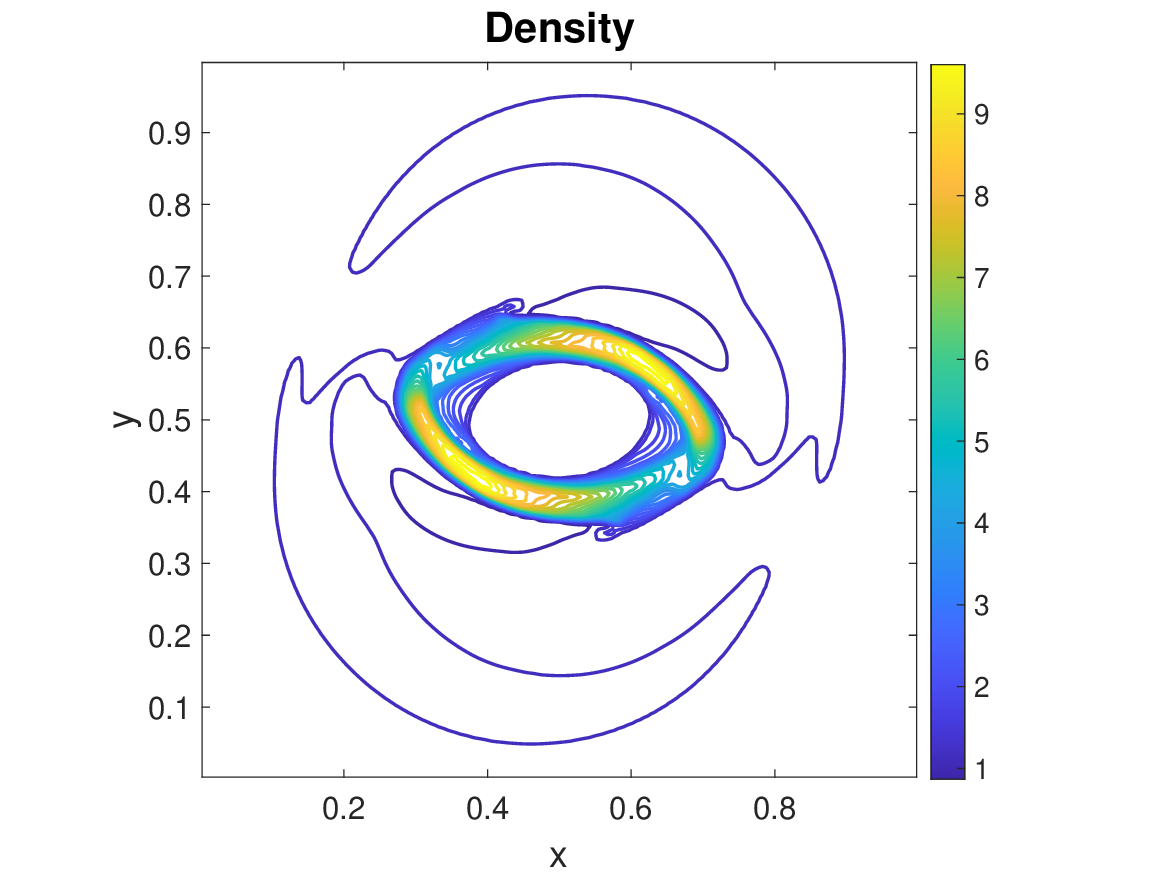}
    }
    \\
    \subfigure[cut along $x=0.5$]{
    \includegraphics[width=0.8\linewidth]{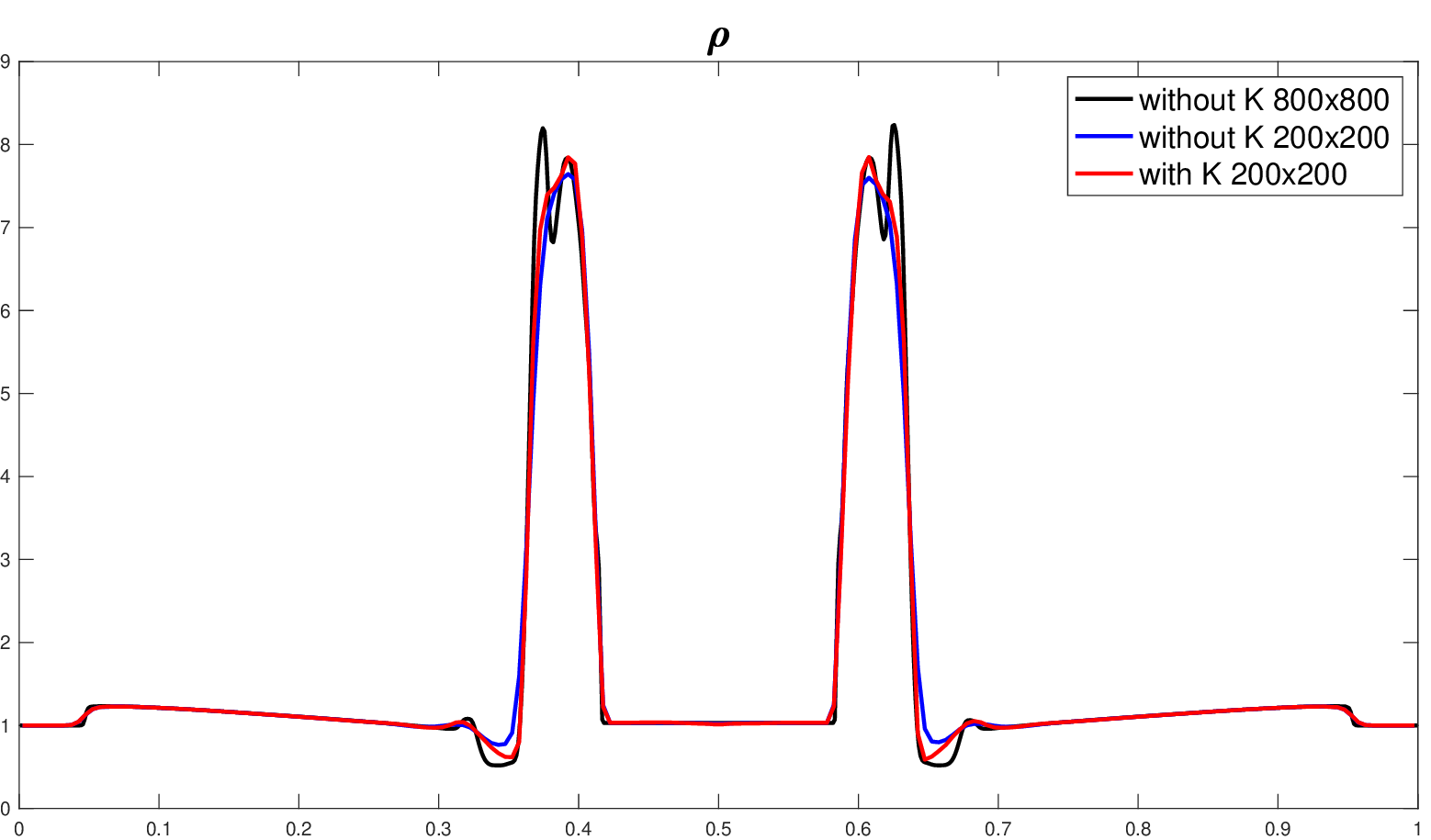}
    }
    \caption{Example~\ref{sec.2Drotor} Rotor test: (a) and (b) are the density with 30 contour lines. (c) is the  comparison of the cross-section along $x=0.5$.}
    \label{fig.rotor_com}
\end{figure}

\begin{figure}[htbp]
    \centering
    \subfigure[Density on $400\times400$]{
        \includegraphics[width=0.48\linewidth]{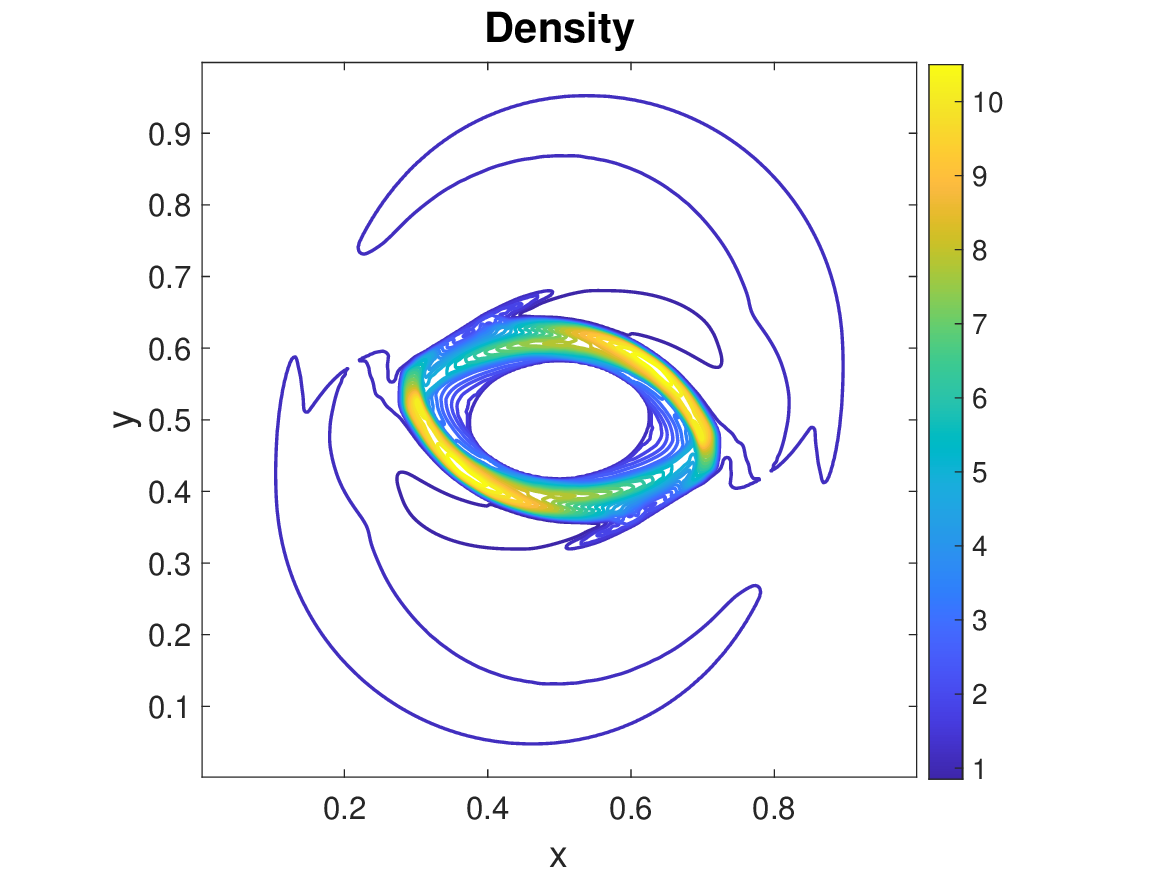}
        } 
    \subfigure[Mach number on $400\times400$]{
        \includegraphics[width=0.48\linewidth]{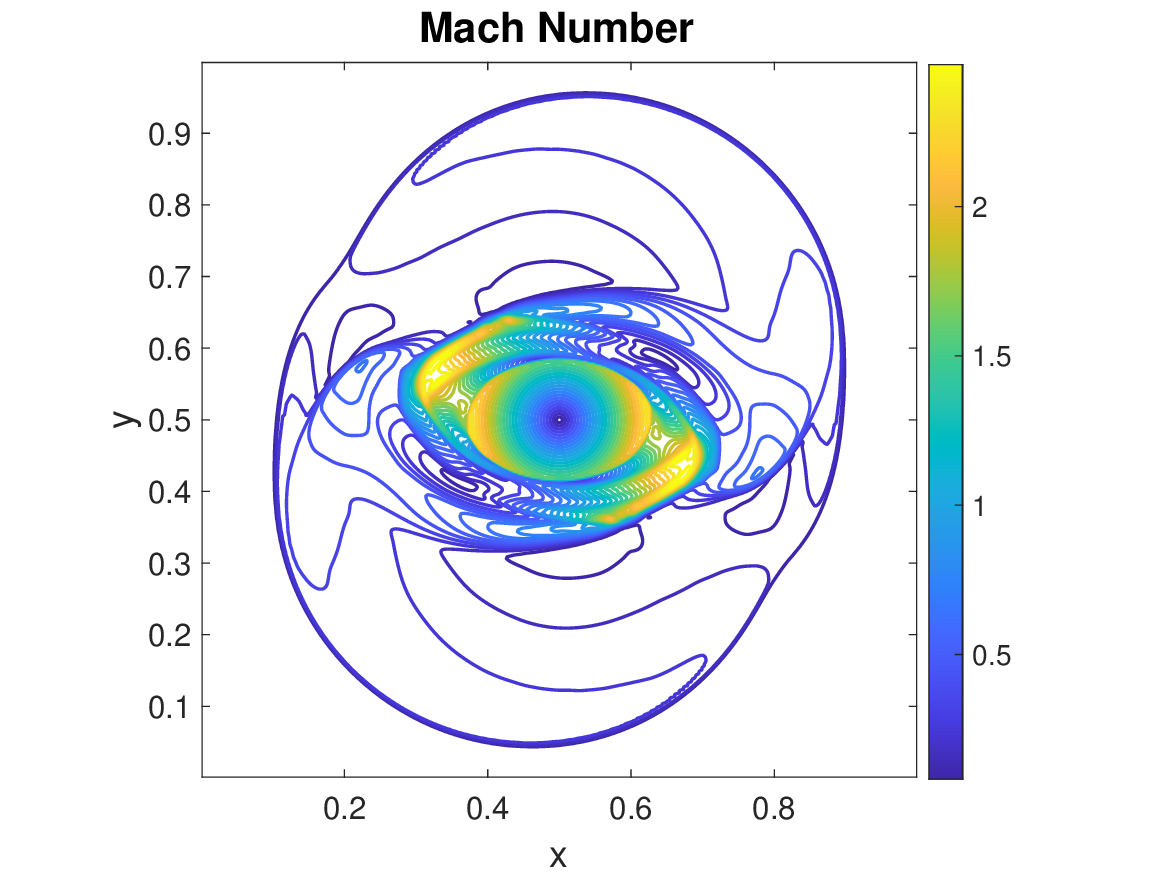}
    }      
    \caption{Example~\ref{sec.2Drotor} Rotor test: Results of density and mach number with 30 contour lines.}
    \label{fig.rotor400}
\end{figure}

The solution develops a rapidly rotating dense fluid disk surrounded by stationary fluid. We show a comparison of the density using the scheme with the correction term K and without K along $x=0.5$ in Figure \ref{fig.rotor_com}. It seems that the result with the correction term performs better, and is closer to the reference solution.  In Figure \ref{fig.rotor400}, we provide the result of density and Mach number on $400\times400$ grid points, which also shows the central rotating part is captured with good circular symmetry. 

\subsubsection{Blast wave}
\label{sec.2DBlast}
\begin{figure}[htbp]
    \centering
    \subfigure[Density]{
        \includegraphics[width=0.48\linewidth]{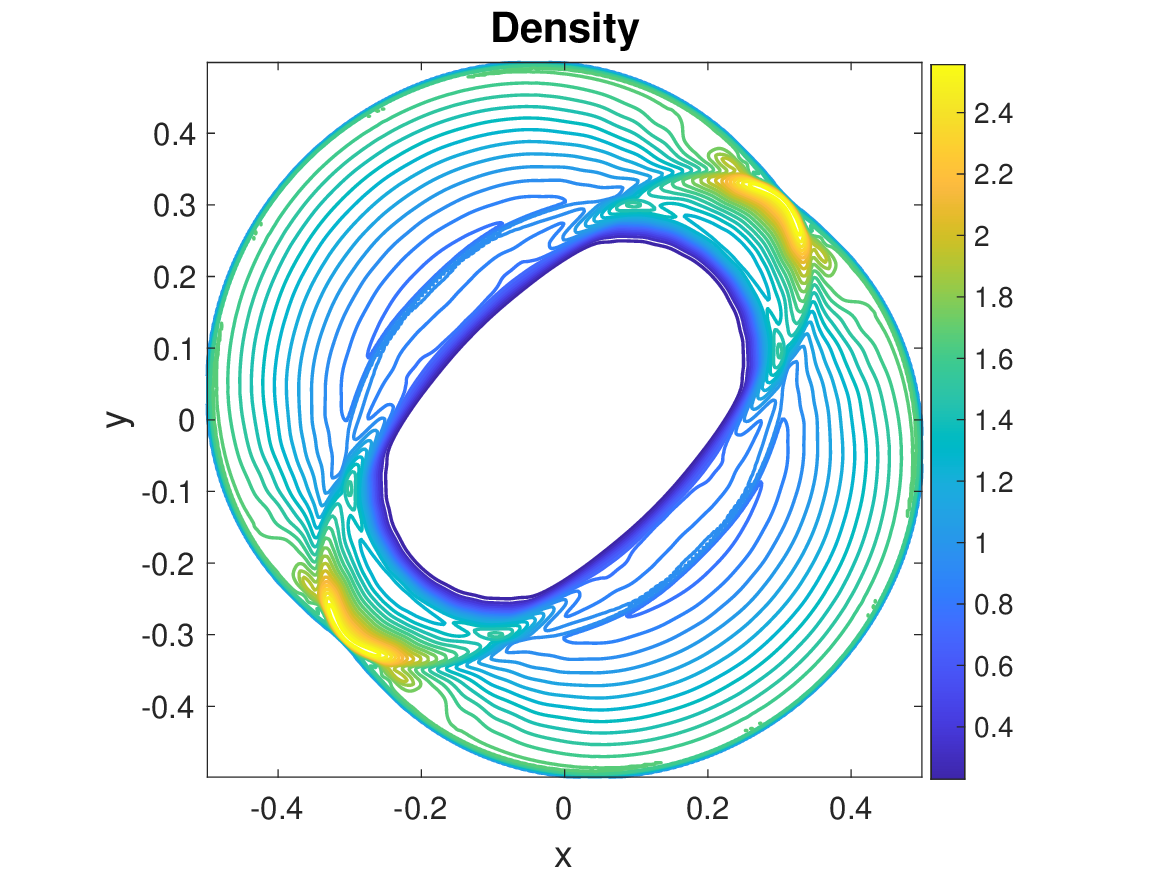}
    } 
    \subfigure[Pressure]{
        \includegraphics[width=0.48\linewidth]{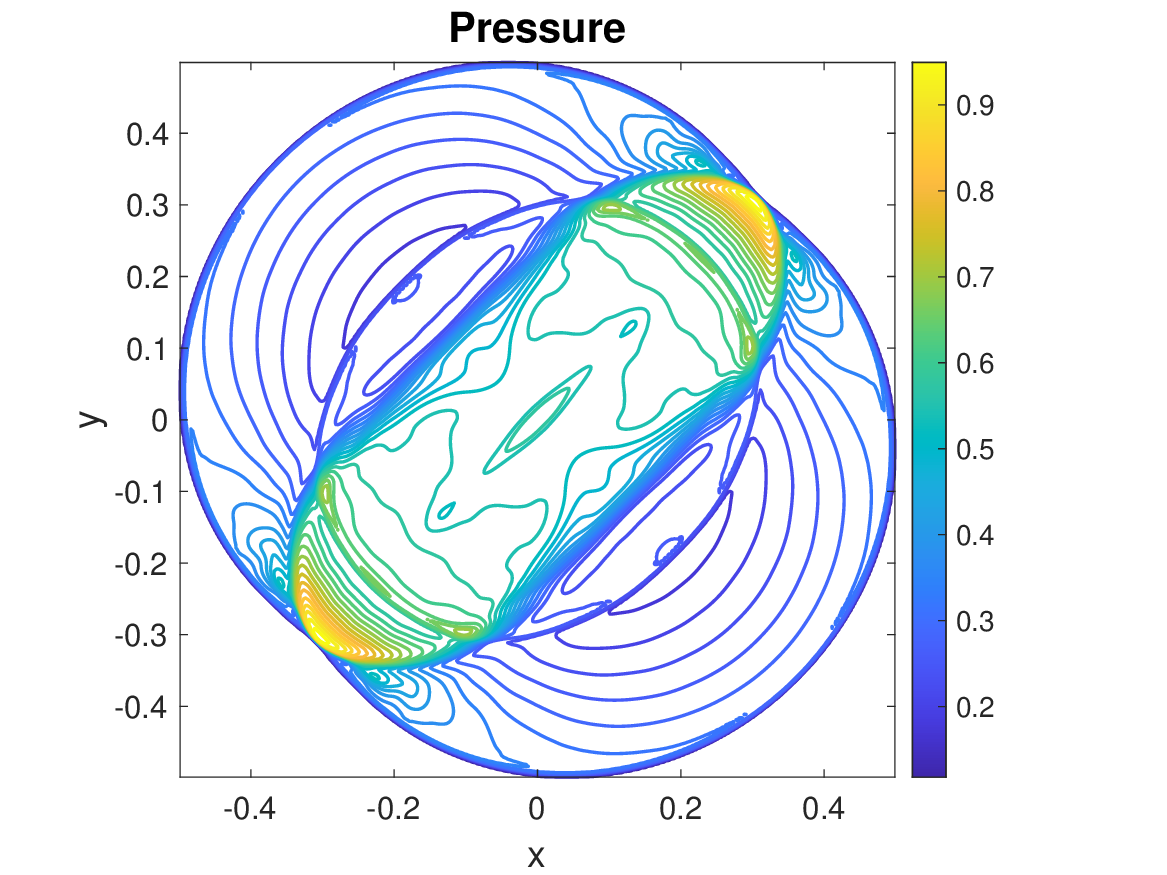}
    } \\
    \caption{Example \ref{sec.2DBlast} Blast wave: Results on $400\times400$ drawn with 30 contour lines.}
    \label{fig.Blast_DP}
\end{figure}

\begin{figure}[htbp]
    \centering
    \subfigure[$400\times400$ grid]{
        \includegraphics[width=0.48\linewidth]{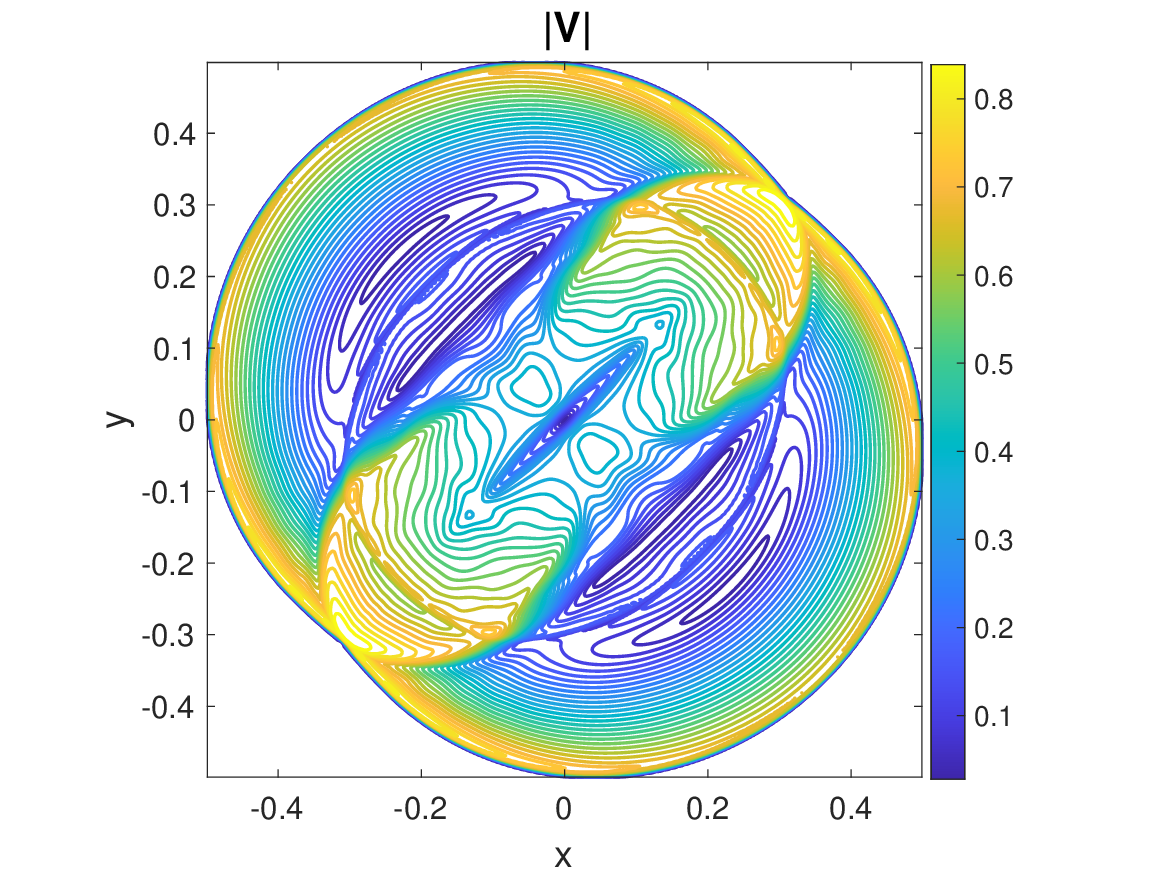}
    } 
    \subfigure[$800\times800$ grid]{
        \includegraphics[width=0.48\linewidth]{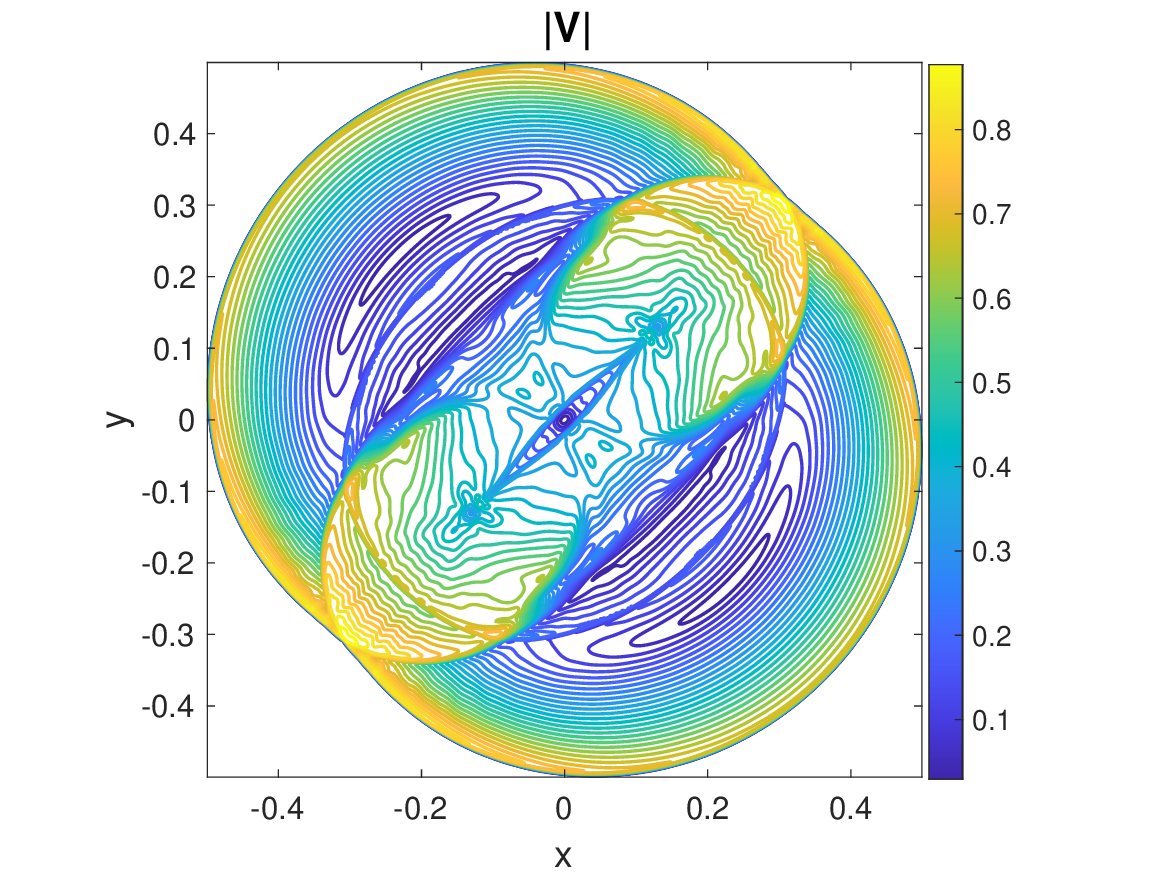}
    } 
    
    \caption{Eaxmple~\ref{sec.2DBlast} Blast wave: Velocity of on $400\times400$ and $800\times800$ with 30 contour lines.}
    \label{fig.Blast_V}
\end{figure}

This problem consists of an initial central region of high pressure which rise to outward going shocks that model a blast wave. We follow the setup in \cite{ref_Blast} and take the computational domain to be $[-0.5,0.5]\times[-0.5,0.5]$ and the final time is $T=0.2$. The initial condition is given by
\begin{equation*}
    (\rho, v_1, v_2, v_3, B_1, B_2, B_3, p)^\top = \left\{
    \begin{array}{cc}
       (1, 0, 0, 0, \dfrac{1}{\sqrt{2}}, \dfrac{1}{\sqrt{2}}, 0, 10)^\top\quad  &\text{for}\;\; \sqrt{x^2+y^2}< 0.1,\\
       (1, 0, 0, 0, \dfrac{1}{\sqrt{2}}, \dfrac{1}{\sqrt{2}}, 0, 0.1)^\top  \;\; &\text{else}.
    \end{array}
    \right. 
\end{equation*}
Figure~\ref{fig.Blast_DP} shows the density and pressure on $400\times400$ grid points, and Figure \ref{fig.Blast_V} shows the velocity on $400\times400$ and $800\times800$ grid points. We can see that our results match the solution in reference \cite{ref_Blast}.

\subsubsection{Challenging Blast problem}
\label{sec.2DCHBlast}
\begin{figure}[htbp]
    \centering
    \subfigure[Density]{
        \includegraphics[width=0.48\linewidth]{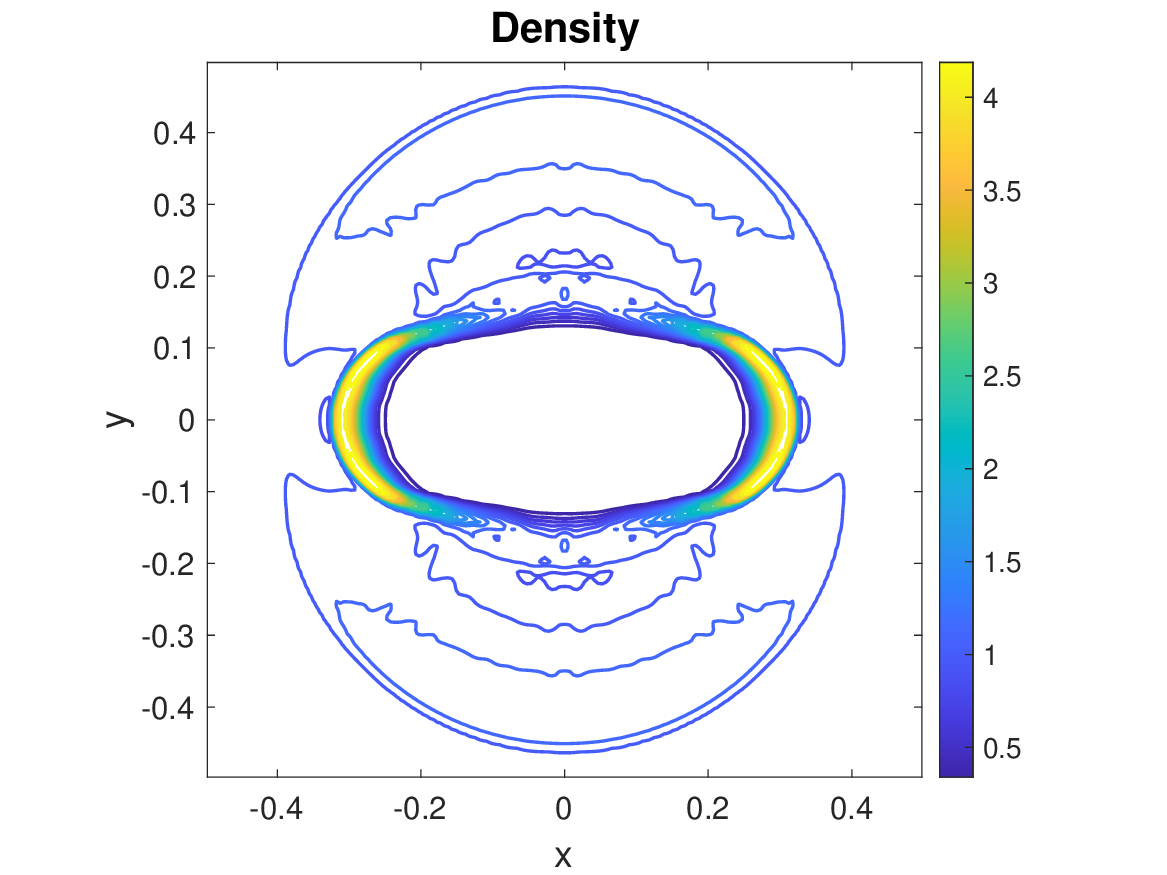}
        } 
    \subfigure[Pressure]{
        \includegraphics[width=0.48\linewidth]{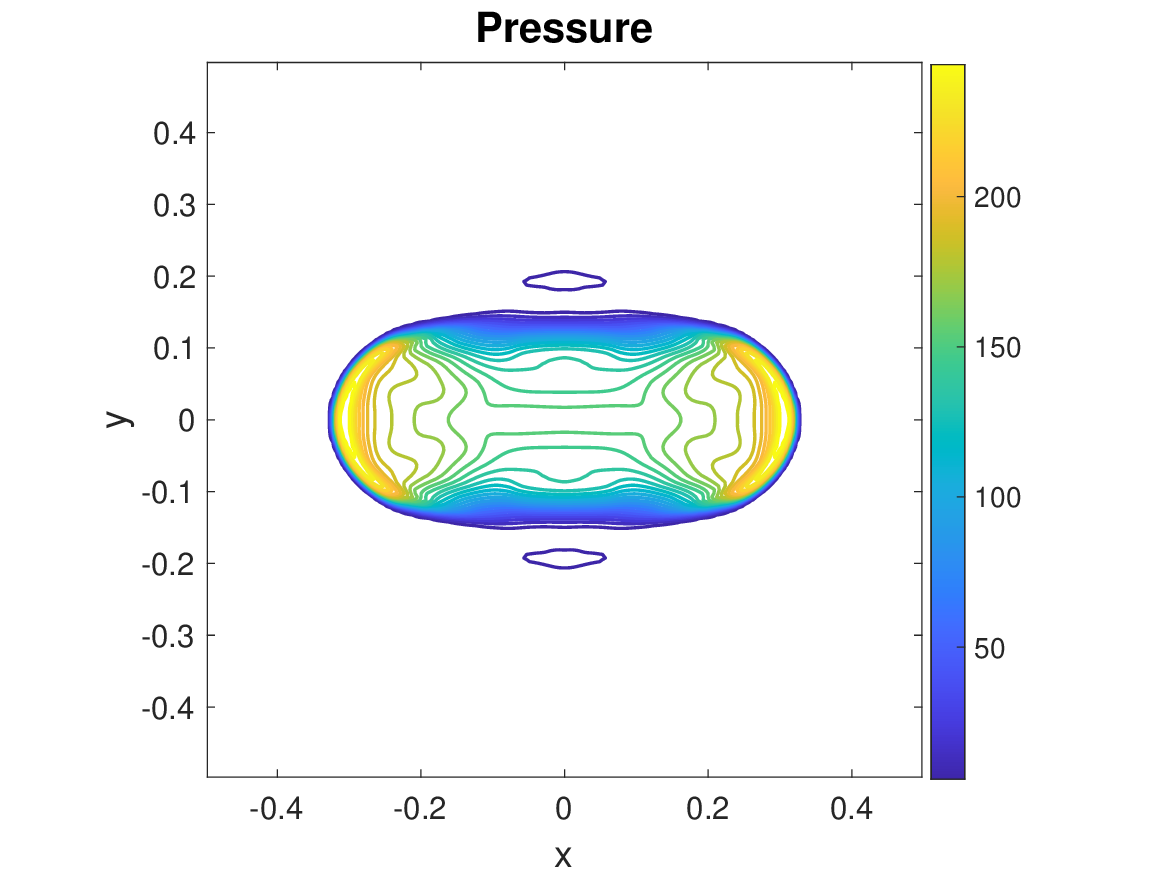}
        } \\
    \subfigure[Velocity]{
        \includegraphics[width=0.48\linewidth]{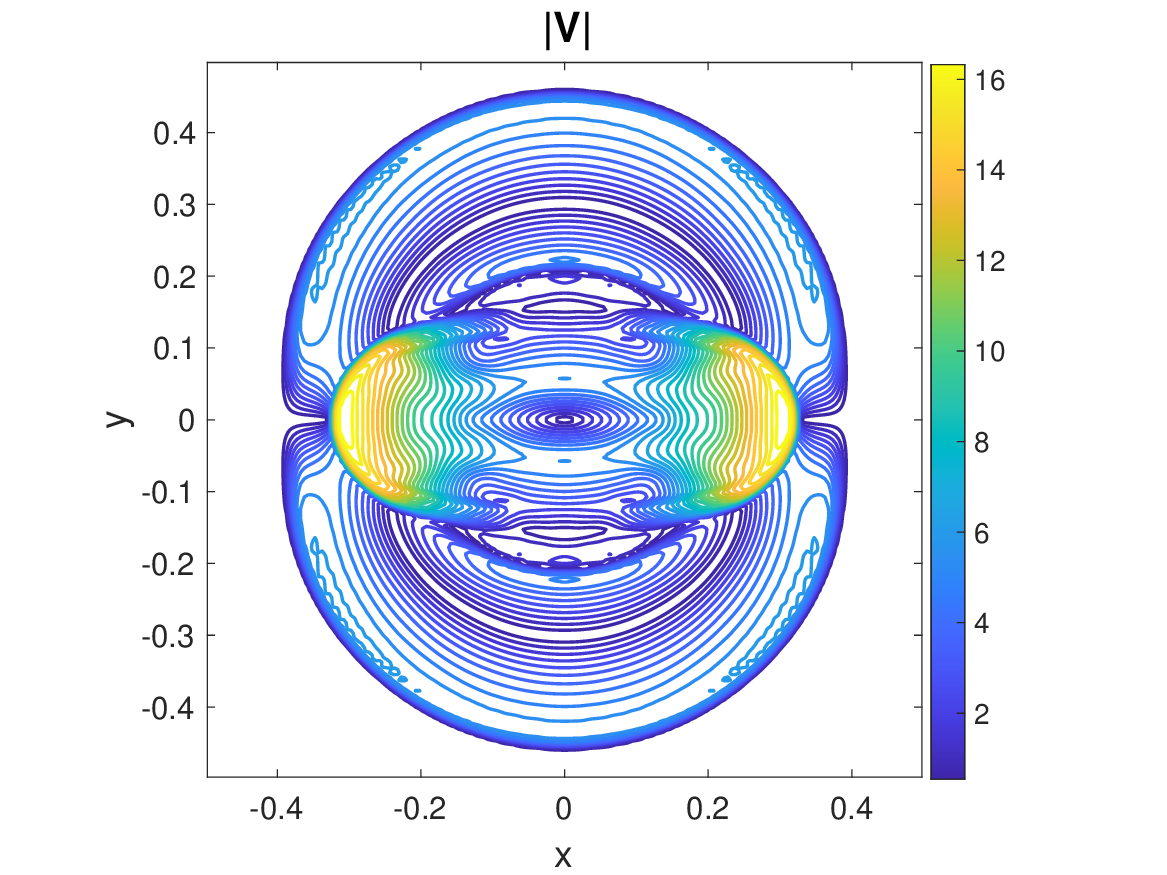}
        } 
    \subfigure[Magnetic pressure]{
        \includegraphics[width=0.48\linewidth]{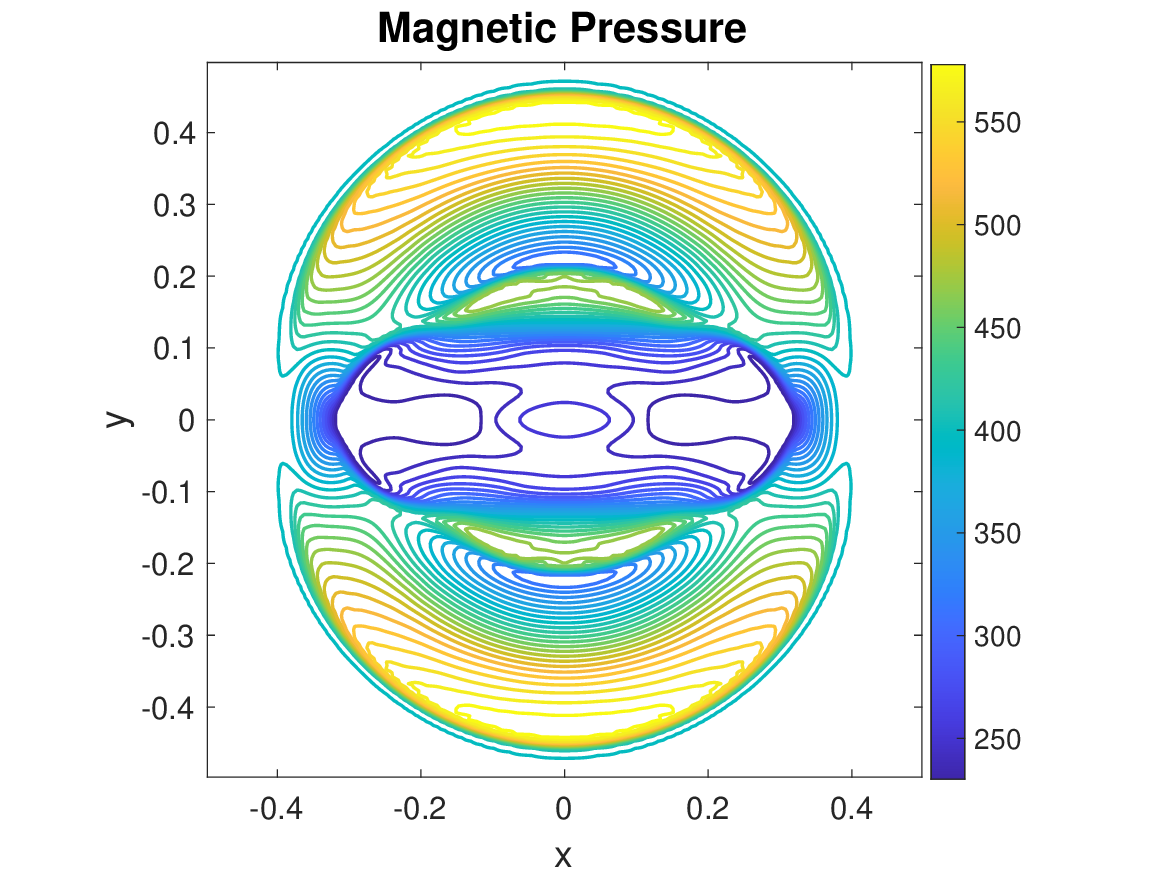}
        } 
    \caption{Example~\ref{sec.2DCHBlast} Challenging Blast problem: Results on $200\times200$ with 30 contour lines.}
    \label{fig.ChBlast}
\end{figure}

This test case is similar to the previous one but gives rise to very low gas pressure regions and strong magnetosonic shocks. This makes it a challenge for numerical methods in terms of preserving positivity of solution. We use the setup in \cite{ref_CHBlast} and compute in the domain $[-0.5, 0.5]\times[-0.5, 0.5]$ with periodic boundary condition. The initial data is
\begin{equation*}
    (\rho, v_1, v_2, v_3, B_1, B_2, B_3, p)^\top = \left\{
    \begin{array}{cc}
       (1, 0, 0, 0, \dfrac{100}{\sqrt{4\pi}}, 0, 0, 1000)^\top\quad  &\text{for}\;\; \sqrt{x^2+y^2}\leq 0.1,\\
       (1, 0, 0, 0, \dfrac{100}{\sqrt{4\pi}}, 0, 0, 0.1)^\top  \;\; &\text{else}.
    \end{array}
    \right. 
\end{equation*}
with $\gamma=1.4$. Figure \ref{fig.ChBlast} shows the results of density, pressure, velocity, and magnetic pressure on the $200\times200$ grid points at the final time $T = 0.01$. The computations remain stable without any loss of positivity and no special treatment has been done to ensure the positivity property. The results are comparable with those in~\cite{ref_CHBlast}.

\section{Summary and conclusions}
In this paper, we apply the concept of the LDCU scheme which was first developed for the Euler system to the MHD system.  We first construct the 1-D fully-discrete MHD scheme, and then develop the semi-discrete MHD scheme. Based on the same idea, we extend our 1-D MHD scheme to the hydrodynamic variables in the 2-D MHD system, and then adopt the constrained transport method for the magnetic variables. The hydrodynamic and magnetic variables are stored in a staggered manner as is standard in all constrained transport methods. With the help of the three-stage Runge-Kutta method, we demonstrate our proposed scheme in several numerical experiments.  The obtained results show that the LDCU correction term in our scheme gives better resolution of contact discontinuity, and can achieve the second order accuracy for smooth solutions. The constrained transport idea helps to keep the divergence close to machine zero which allows stable computations for challenging problems including the blast problems without issue of positivity violation. These ideas lead to a simple, Riemann solver-free method that is able to compute MHD problems in a stable and accurate manner.

\section*{Acknowledgments}
The authors are grateful to Dr. Junming Duan for his constructive feedback and expert advice on the manuscript. The work of Praveen Chandrashekar is supported by the Department of Atomic Energy, Government of India, under project no.~12-R\&D-TFR-5.01-0520.

\bibliographystyle{alpha}
\bibliography{reference}

\begin{appendices}
\section{Algorithm}\label{sec:algo}
\begin{algorithm}
\caption{2-D LDCU MHD scheme}
\begin{algorithmic}[1]
\Require{Initial variables $\rho, v_1, v_2, v_3, B_1, B_2, B_3, p$ at cell center, and $B_{1,\bf{x}_f}, B_{2,\bf{y}_f}$ at interface} 
\Ensure{Variables at final time $T_f$}
\Statex
\State Set ghost cells around the boundaries
\For{each face}
\State Reconstruct variables $\rho, \rho v_1, \rho v_2, \rho v_3, B_1, B_2, B_3, E$ with limiter
\State Compute the maximal local wave speed $a^{\pm},b^{\pm}$
\EndFor
\State Compute time step $\Delta t$ and set time counter $t=\Delta t$
 \While{$t\leq T_f$}
  \For{each RK stage}
  \For{each face}
  \State Compute the LDCU term $K^*$
  \State Compute flux $F$, and $G$
  \State Compute the average velocities $\overline{V}_{2,\bf{x}_f}$ and $\overline{V}_{1,\bf{y}_f}$ at interface
  \EndFor
  \For{each vertices}
  \State Reconstruct $\overline{V}_{2,\bf{x}_f}, \overline{V}_{1,\bf{y}_f}$, and $B_{1, \bf{x}_f}, B_{2,\bf{y}_f}$ to the corner
  \State Compute electric field $\Omega$ 
  \EndFor
  \For{each face}
  \State Update variable $B_{1,\bf{x}_f}$ and $B_{2,\bf{y}_f}$ to next RK stage
  \EndFor
  \For{each cell}
  \State Update variable $\rho, \rho v_1, \rho v_2, \rho v_3, B_3, E$ to next RK stage
  \State Compute $B_1$ and $B_2$ at next RK stage
  \EndFor
  \State Repeat line 1 to 5
  \EndFor
  \State Compute time step $\Delta t$ and do $t=t+\Delta t$  
  \EndWhile
\end{algorithmic}
\end{algorithm}

\end{appendices}

\end{document}